\documentclass[12pt]{amsart}
\usepackage{amsmath}
\usepackage{amssymb}
\newtheorem{theorem}{Theorem}[section]
\newtheorem{lemma}[theorem]{Lemma}
\newtheorem{proposition}[theorem]{Proposition}
\newtheorem{corollary}[theorem]{Corollary}
\newtheorem{addendum}[theorem]{Addendum}

\theoremstyle{definition}
\newtheorem{definition}[theorem]{Definition}
\newtheorem{definitions}[theorem]{Definitions}
\newtheorem{example}[theorem]{Example}

\newtheorem{definitions and remarks}[theorem]{Definitions and Remarks}

\theoremstyle{remark}
\newtheorem{remark}[theorem]{Remark}
\newtheorem{remarks}[theorem]{Remarks}
\newtheorem{questions}[theorem]{Questions}
\newtheorem{algorithm}[theorem]{Algorithm}

\numberwithin{equation}{section}

\newcommand{\inv}{\mathrm{inv}}
\newcommand{\rk}{\mathrm{rk}\,}
\newcommand{\Sing}{\mathrm{Sing}\,}
\newcommand{\supp}{\mathrm{supp}\,}
\newcommand{\mon}{\mathrm{mon}\,}

\newcommand{\codim}{\mathrm{codim}\,}
\newcommand{\car}{\mathrm{char}\,}

\newcommand{\length}{\mathrm{length}\,}

\newcommand{\al}{{\alpha}}
\newcommand{\be}{{\beta}}
\newcommand{\de}{{\delta}}
\newcommand{\ep}{{\epsilon}}
\newcommand{\D}{{\Delta}}
\newcommand{\g}{{\gamma}}
\newcommand{\la}{{\lambda}}
\newcommand{\La}{{\Lambda}}
\newcommand{\Om}{{\Omega}}

\newcommand{\s}{{\sigma}}
\newcommand{\Sig}{{\Sigma}}
\newcommand{\vp}{{\varphi}}
\newcommand{\io}{{\iota}}
\newcommand{\om}{{\omega}}
\newcommand{\Th}{{\Theta}}

\newcommand{\IN}{{\mathbb N}}

\newcommand{\IA}{{\mathbb A}}
\newcommand{\IR}{{\mathbb R}}

\newcommand{\IT}{{\mathbb T}}
\newcommand{\IZ}{{\mathbb Z}}
\newcommand{\IF}{{\mathbb F}}

\newcommand{\cA}{{\mathcal A}}
\newcommand{\cC}{{\mathcal C}}
\newcommand{\cD}{{\mathcal D}}
\newcommand{\cE}{{\mathcal E}}

\newcommand{\cG}{{\mathcal G}}
\newcommand{\cH}{{\mathcal H}}
\newcommand{\cI}{{\mathcal I}}

\newcommand{\cO}{{\mathcal O}}
\newcommand{\cQ}{{\mathcal Q}}

\newcommand{\uk}{\underline{k}}
\newcommand{\um}{{\mathfrak m}}
\newcommand{\of}{\overline{f}}

\newcommand{\tal}{{\tilde \al}}
\newcommand{\tbe}{{\tilde \be}}
\newcommand{\tg}{{\tilde \g}}

\newcommand{\wcO}{{\widehat \cO}}
\newcommand{\wR}{{\widehat R}}
\newcommand{\wI}{{\widehat I}}

\newcommand{\ucG}{\underline{\cG}}
\newcommand{\ucH}{\underline{\cH}}
\newcommand{\ucC}{\underline{\cC}}
\newcommand{\ucD}{\underline{\cD}}
\newcommand{\ucQ}{\underline{\cQ}}
\newcommand{\uzero}{\underline{0}}
\newcommand{\uone}{\underline{1}}

\newcommand{\frN}{{\mathfrak N}}
\newcommand{\frV}{{\mathfrak V}}

\newcommand{\lbr}{{[\![}}
\newcommand{\rbr}{{]\!]}}

\begin{document}

\title[Desingularization of toric and binomial varieties]
{Desingularization of\\ 
toric and binomial varieties}

\author{Edward Bierstone}
\address{Department of Mathematics, University of Toronto, Toronto,
Ontario, Canada M5S 3G3}
\email{bierston@math.toronto.edu}
\thanks{The authors' research was supported in part by NSERC 
grants OGP0009070 and OGP0008949.}

\author{Pierre D. Milman}
\address{Department of Mathematics, University of Toronto, Toronto,
Ontario, Canada M5S 3G3}
\email{milman@math.toronto.edu}



\begin{abstract}
We give a combinatorial algorithm for equivariant embedded resolution of 
singularities of a toric variety defined over a perfect field.
The algorithm is realized by a finite succession of blowings-up 
with smooth invariant centres that satisfy the normal flatness condition
of Hironaka. The results extend to more general varieties defined 
locally by binomial equations.
\end{abstract}

\maketitle
\setcounter{tocdepth}{1}
\tableofcontents

\section{Introduction}
We give a simple combinatorial algorithm for equivariant
embedded resolution of singularities of a toric variety $X$
(not necessarily normal) defined over a perfect field $\uk$.
The algorithm is realized by a finite
succession of blowings-up with smooth invariant centres, 
such that each successive transform 
of $X$ is normally flat along the corresponding centre (condition 
of Hironaka \cite{Hann}).  We
announced this article in \cite{BMda1} as ``Desingularization 
algorithms II. Toric and locally binomial varieties'', but have
changed the title because it is largely independent of
\cite{BMda1}.

Throughout this paper, $\uk$ denotes a field. If $X$ is a
toric variety over $\uk$, we denote by $T_X$ the embedded algebraic
torus. We say that
a morphism of toric varieties $f: X \to Y$ is {\it equivariant} 
if $f|T_X$ is a homomorphism of tori $h: T_X \to T_Y$ and $f$ is
equivariant with respect to the homomorphism $h$.
(Some basic facts about toric varieties are recalled
in Section 2 below.)

\begin{theorem}[{\rm Equivariant embedded desingularization of a 
toric variety}] 
Let $X \hookrightarrow M$ denote an equivariant
embedding of toric varieties over a perfect field $\uk$, 
where $M$ is smooth.
Then there is a finite sequence
of blowings-up of $M$,
\begin{equation}
M = M_0 \stackrel{\pi_1}{\longleftarrow} M_1 \longleftarrow \cdots
\stackrel{\pi_{t+1}}{\longleftarrow} M_{t+1}\ ,
\end{equation}
such that:
\begin{enumerate}
\item 
The centre $D_j$ of each blowing-up $\pi_{j+1}$ is a smooth
$T_{M_j}$-invariant subvariety of $M_j$.
\item
Set $X_0 = X$. For each $j = 0,\ldots,t$, let $X_{j+1}$ denote
the strict transform of $X_j$ by $\pi_{j+1}$.
Then each $C_j := D_j \cap X_j$ is a smooth $T_{X_j}$-invariant 
subvariety
of $X_j$, and $X_j$ is normally flat along $C_j$.
\item
For each $j = 0,\ldots,t$, let $E_{j+1}$ denote the exceptional
divisor of $\pi_1\circ\cdots\circ\pi_{j+1}$. Then, for each
$j = 0,\ldots,t$, either 
$C_j \subset \Sing X_j$, or $X_j$ is smooth
and $C_j \subset X_j \cap E_j$.
\item
$X_{t+1}$ is smooth, and $X_{t+1}$, $E_{t+1}$ simultaneously have 
only normal crossings.
\end{enumerate}
\end{theorem}

The condition $C_j = D_j \cap X_j$ (intersection as subspaces, or
subschemes) implies that $X_{j+1} =$ blowing-up of $X_j$ with
centre $C_j$ \cite[\S II.7]{Hart}. {\it Normal flatness} of $X_j$ 
along $C_j$ means the normal cone of $C_j$ in $X_j$ is flat over 
$C_j$; this is 
equivalent to the property that the 
{\it Hilbert-Samuel function} $H_{X_j,x}: \IN \to \IN$,
$$
H_{X_j,x}(l)\ :=\ \length \frac{\cO_{X_j,x}}{\um_{X_j,x}^{l+1}}\ ,
$$
is locally constant as a function of $x \in C_j$ \cite{Be}.
($\um_{X_j,x}$ denotes the maximal ideal of $\cO_{X_j,x}$.)

A sequence of blowings-up as in the
theorem has the property that, for each $j$, $T_{M_j} = T_M$,
$T_{X_j} = T_X$,
and the embedding $X_j \hookrightarrow M_j$
is equivariant. 

For each $j$, every irreducible component of $C_j$ is a smooth
orbit closure of $X_j$. We can assume that each $D_j$ is the union
(necessarily disjoint) of the smallest orbit closures of $M_j$
containing the components of $C_j$. Although Theorem 1.1 could be
restated with each $C_j, D_j$ taken to be orbit closures, we need
to allow disconnected centres to have also the following.

\begin{addendum}[{\rm Canonical equivariant embedded desingularization}]
For every equivariant embedding $X \hookrightarrow M$ of toric varieties
over a perfect field $\uk$, 
where $M$ is smooth, there is sequence of blowings-up (1.1)
satisfying the conditions of Theorem 1.1, with the following property:
If $\iota: M' \to M$ is an open equivariant embedding, then 
$X$ and $X' = \iota^{-1}(X)$ have the same resolution towers over $M'$
(not counting isomorphisms in the sequences of blowings-up).
\end{addendum}

We prove Theorem 1.1 in Section 8 below. Desingularization can be
realized by a simple combinatorial algorithm that is canonical
up to an ordering of the codimension $1$ orbit closures of $M$ (\S9.1).
(Every orbit closure of a smooth toric variety is an intersection of
codimension $1$ orbit closures.)
The canonical desingularization algorithm of Addendum 1.2 comes
at an extra cost: Additional blowings-up are needed to replace
the codimension $1$ orbit closures involved 
by components of exceptional divisors 
(which are
ordered by the sequence of blowings-up). (See \S9.2.)

At a comparable cost, we can give a canonical desingularization 
algorithm for more general {\it toroidal} or
{\it binomial} varieties defined over $\uk$. For these classes,
the locally toric or binomial structures are related by global
divisors (Section 10). 

The general desingularization algorithm
of \cite{BMinv, BMda1} can be adapted to varieties over
$\uk$ that are merely locally toric or locally binomial, but again
has additional complexity due to blowings-up needed to guarantee that
combinatorial centres chosen locally as in the
toric case extend to global smooth centres. The general algorithm,
however, strengthens condition (3) in Theorem 1.1: $C_j$ and $E_j$
will simultaneously have only normal crossings, for all $j$. (The
simpler toric or binomial algorithms provide only simultaneous normal
crossings of $C_j$ and $N \cap E_j$, where $N$ is a local minimal
embedding variety of $X_j$.)

\begin{example}[{\rm cf. \cite[Example 1.2]{BMda1}}]
Consider the toric hypersurface $X$ in $\uk^3$
defined by
$$
z^d - x^{d-1}y^d \ = \ 0\ .
$$
Let $n($T$)$, $n($CT$)$, $n($CB$)$
and $n($BM$)$ denote the number of blowings-up needed to reduce the
maximum order $d$ of $X$, using our toric desingularization
algorithm, canonical toric algorithm, canonical binomial algorithm, 
and the general desingularization
algorithm of \cite{BMinv, BMda1}, respectively. Then
\begin{align*}
n(\mbox{T})\ &=\ 1\\
n(\mbox{CT}),\ n(\mbox{CB})\ &\leq\ d + j\\
n(\mbox{BM})\ &\leq\ 2d + k \, ,
\end{align*}
where $j$ and $k$ are independent of $d$.
\end{example}

The generalizations of our results on toric varieties are presented
briefly in Section 10, but Section 6 and \S\S7.1-7.3 already
cover the more general binomial ideals, to prepare the ground.

An application of our results to desingularization of the induced
metric and of the Gauss mapping of an embedded binomial variety
will be presented in a forthcoming paper.

Our toric desingularization algorithm can be viewed as the
combinatorial part of the general algorithm of \cite{BMinv, BMda1}.
From the point of view of resolution of singularities over
fields of arbitrary characteristic, 
the classes of toric or locally binomial varieties
present a great simplification because of the trivial existence
of smooth subvarieties of {\it maximal contact}; cf. Example 2.5 below.
Our techniques are of a ``differential calculus'' nature. The
hypothesis that $\uk$ is perfect is essentially equivalent to
the possibility of using calculus in local coordinates.
(See Section 3.) Of course, every algebraically closed field
is perfect.

Combinatorial desingularization of normal toric varieties goes
back to \cite{KKMS}. A normal toric variety $X$ corresponds to a 
fan $\Sig$ in a lattice; $X = X(\Sig)$ is smooth if and only if
$\Sig$ is a regular fan (cf. \S2.2 below). Any fan can be
refined to a regular fan by iterated star-subdivisions
\cite[\S2.6]{Ful}, \cite[\S5]{Cox}; 
the corresponding morphism of toric varieties is a
resolution of singularities (i.e., a proper morphism from a smooth
variety to $X$ that is an isomorphism over the complement of the
singular locus). Each star-subdivision corresponds to a normalized
blowing-up (with centre not necessarily smooth). In the case of an 
equivariant embedding
$X \hookrightarrow M$, this provides, in general, 
neither equivariant
embedded desingularization, nor is the morphism given by 
(or evidently dominated by) a sequence of
blowings-up with centers satisfying the conditions of our theorem. 
(See \cite[\S6.2]{Tei}.) De Concini and Procesi have shown
that if $M'$ is a smooth toric variety birationally equivalent to
$M$, then there is an equivariant birational morphism $M'' \to M'$,
where $M''$ is obtained from $M$ by a finite sequence of blowings-up 
with centres that are codimension $2$ orbit closures \cite{DP}.

Theorem 1.1 (and Addendum 1.2) can be simplified in the case
that $X \hookrightarrow M$ is a toric {\it hypersurface} in $M$
(i.e., $X$ is defined by a principle ideal $\cI_X$ in $\cO_M$).
Moreover, in this case, we do not need to assume that $\uk$ is perfect.

\begin{theorem}[{\rm Equivariant embedded desingularization of a 
toric hypersurface}] 
Let $X \hookrightarrow M$ denote an equivariant
embedding of toric varieties over a field $\uk$, where $M$ is smooth.
Assume that $X$ is a hypersurface in $M$. Then Theorem 1.1
(and Addendum 1.2) can be strengthened so that each $\pi_{j+1}$ 
is a blowing-up with $T_M$-invariant centre $C_j = D_j$ in $X_j$.
\end{theorem}

In the case of a hypersurface, the condition that $X_j$ is 
normally flat along $C_j$
is equivalent to the condition that the order of
vanishing of $\cI_{X_j}$ is locally constant on $C_j$.

A smooth toric variety $M$ corresponds to a regular fan $\Sig$ in 
a lattice $L \cong \IZ^n$. A blowing-up of $M$ with smooth
$T_M$-invariant centre correponds to a certain star-subdivision
of $\Sig$. (See Section 4.)
A toric hypersurface $X \hookrightarrow M = M(\Sig)$ is
determined by (the restriction to $\Sig$ of) an integer-valued
linear function $\la$ on $M$. (See Lemma 2.6.) Theorem 1.4 translates
into a purely combinatorial statement about linear extension of
$\la$ to successive star-subdivisions of $\Sig$ (Theorem 5.3.)

Our proof of Theorem 1.4 shows that (if $X_j$ is not already 
smooth) the locus of maximum order of $\cI_{X_j}$ has only
normal crossings and, morevover, each of its irreducible components 
is an orbit-closure of $M$. We can resolve the singularities by
choosing any of these components as the centre $C_j$ of the
next blowing-up. (Compare with Theorem 1.6.)

A general toric subvariety $X \subset M$ 
cannot be desingularized by blowings-up with 
$T_M$-{\it invariant} centres that lie in the successive
strict transforms of $X$, as in the hypersurface case. 
(We say that a toric variety $X \subset M$
is a {\it toric subvariety} of $M$ if the inclusion $X \hookrightarrow
M$ is an equivariant embedding of toric varieties.)

\begin{example}
Let $X \subset \uk^6$ denote the affine toric variety
whose ideal is generated by the binomials
$$
f \ =\ w^2 - uv\ , \quad g\ =\ z - xy
$$
(in six variables $u,v,w,x,y,z$). Let $N \subset \uk^6$ denote the
smooth subvariety $\{g = 0\}$. Then $X \hookrightarrow N$ is a
minimal embedding; $N$ is a toric subvariety of $M$. Desingularization
of $X$ is the result of a single blowing-up with centre $C$ defined by
$$
u\ =\ v\ =\ w\ =\ 0\ , \quad z - xy\ =\ 0\ .
$$
The centre $C$ is $T_X$- or $T_N$-invariant, 
but is not $(\uk^*)^6$-invariant.
(It is not a coordinate subspace of $\uk^6$.) Let $D = \{u = v = w =
0 \} \subset \uk^6$. Then $D$ is the smallest orbit closure of 
$\uk^6$ containing $C$; $D$ and $N$ have only normal crossings,
and $D \cap N = C$.
\end{example}

An affine toric subvariety of $\uk^n$ corresponds to a prime ideal
in $\uk[x_1,\ldots,x_n]$ generated by binomials (a {\it toric
ideal}; see Section 2 below). A toric ideal has a distinguished
set of binomial generators (that we call a {\it standard basis}),
uniquely determined by an ordering of the variables (Theorem 6.2).
A standard basis is similar to the notion of Gr\"obner basis,
with the difference that the initial monomial of an element of
a standard basis is a monomial of lowest degree, rather than
of highest degree as in a Gr\"obner basis. 

Standard bases play a key role in resolution of singularities
because, in the general case,
we can work with the elements of a standard basis in
the same way as with a single defining equation in the hypersurface
case (or with an arbitrary system of generators in ``principalization
of an ideal''). The properties needed are made precise in Theorems 7.1,
7.3: The locus of maximum order of the elements of a standard basis 
coincides with the maximal locus of the Hilbert-Samuel function
(the ``maximal Samuel stratum''). Moreover, if the Hilbert-Samuel
function does not decrease on blowing-up with centre in the Samuel
stratum, then the standard basis transforms to the standard basis
of the ideal of the strict transform. The following theorem is a
corollary of the first assertion (proof in \S7.4).

\begin{theorem}
Let $X$ denote a toric variety over a perfect field $\uk$.
Then the maximal Samuel stratum $S$ of $X$ has
only normal crossings. Moreover, each irreducible component of
$S$ is a closed smooth $T_X$-invariant subspace of $X$.
\end{theorem}

Unlike the situation in Theorem 1.4 on the hypersurface case,
however, it is not, in general, possible to resolve singularities
by choosing as each successive centre of blowing up an arbitrary
component of the maximal Samuel stratum (Example 7.6).

We do not know whether one can get an algorithm for equivariant
embedded desingularization of a toric variety by blowing up with
centre given at each step by {\it some} component of the maximal
Samuel stratum. (We doubt that this is true.) Our desingularization
algorithm involves blowing up with perhaps smaller
centres that serve to order monomials appearing in a standard
basis, in preparation for blowing up components of the Samuel
stratum as in the hypersurface case.

\section{Toric varieties}
Let $\uk$ denote a field. A {\it toric variety (over $\uk$)} is
an algebraic variety $X$ over $\uk$ that contains an algebraic torus
$T = T_X$ as an open dense subset, and has an action $T \times
X \to X$ of $T$ that extends the natural action of $T$
on itself. {\it We do not assume that $X$ is necessarily normal.}
(In general, an {\it algebraic variety} over $\uk$ means a reduced
scheme of finite type over $\uk$.)

\subsection{Affine toric varieties}
An affine toric variety is simply an affine variety that is 
parametrized by a set of Laurent monomials. Consider a subset $\cA = 
\{a_1, \ldots, a_n\} \in \IZ^d$. Each vector $a_i = (a_{1i}, \ldots,
a_{di})$ identifies with a monomial $t^{a_i} = t_1^{a_{1i}} \cdots
t_d^{a_{di}}$ in the Laurent polynomial ring $\uk[t^\pm] =
\uk[t_1, \ldots, t_d, t_1^{-1}, \ldots, t_d^{-1}]$. We define
the {\it toric ideal} $I_\cA \subset \uk[x] = \uk[x_1, \ldots, x_n]$
as the kernel of the algebra homomorphism 
$$
\uk[x] \to \uk[t^\pm]\ , \quad
x_i \mapsto t^{a_i}\ .
$$
(See \cite[Chs. 4, 13]{Stu2}.)

Consider the group homomorphism
$$
\pi:\ \IZ^n \to \IZ^d\ , \quad 
\g = (\g_1, \ldots, \g_n) \mapsto \g_1 a_1 + \cdots + \g_n a_n\ .
$$
Any $\g \in \IZ^n$ can be written uniquely as $\g = \g^+ - \g^-$,
where $\g^+$ and $\g^- \in \IZ^n$ are nonnegative and have disjoint
supports. See \cite{Stu2} for the following three lemmas.

\begin{lemma}
The toric ideal $I_\cA$ is spanned (as a $\uk$-vector space)
by the binomials
$$
x^{\g^+} - x^{\g^-}, \quad \g \in \ker \pi\ .
$$
\end{lemma}

(A {\it binomial} means a difference of two monomials.) 

\begin{lemma}
An ideal in $\uk[x]$ is toric if and only if it is prime and
generated by binomials.
\end{lemma}

Let $X_\cA \subset \IA^n$ denote the affine variety $V(I_\cA)$
defined by the toric ideal $I_\cA$ (where $\IA^n$ denotes 
$n$-dimensional affine space over $\uk$). A variety of the form $X_\cA$
is an affine toric variety. (See below.)

Let $\rk \cA$ denote the rank of the $d \times n$ matrix with
columns $a_1, \ldots, a_n$ (so that $\rk \cA$ is the dimension
of the lattice $\IZ \cA \subset \IZ^d$ spanned by $a_1, \ldots, a_n$).

\begin{lemma}
The Krull dimension of $\uk[x]/I_\cA$ equals $\rk \cA$.
\end{lemma}

\begin{proof}
The ring $\uk[x]/I_\cA$ is isomorphic to the subring 
$\uk[t^{a_1}, \ldots, t^{a_n}]$ of $\uk[t^\pm]$. The Krull dimension
of the latter is the maximum number of algebraically independent
monomials $t^{a_i}$. But a set of monomials ${t^{a_i}}$ is
algebraically independent if and only if the exponent vectors
${a_i}$ are linearly independent (by Lemma 2.1).
\end{proof}

\begin{example}
Let $\IT$ denote the algebraic group whose $\uk$-rational points
are the points of the multiplicative group $\uk^*$ of non-zero 
elements of $\uk$. The algebraic torus $\IT^n$ has the structure 
of an affine subvariety $V(I)$ of $\IA^{2n}$
determined by the ideal $I \subset \uk[y,z] = 
\uk[y_1, \ldots, y_n, z_1,
\ldots, z_n]$ generated by $y_j z_j - 1$, $j = 1, \ldots, n$.
($\uk[y,z]/I \cong \uk[y^\pm]$.)
\end{example}

An {\it algebraic torus} is an algebraic group isomorphic to 
$\IT^n$, for some $n\geq 1$.

Consider the action of the torus $\IT^d$ on $\IA^n$ (or $(\uk^*)^d$
on $\uk^n$) given by
$$
(t, x) = (t_1, \ldots, t_d, x_1, \ldots, x_n) \mapsto 
(t^{a_1}x_1, \ldots, t^{a_n}x_n)\ .
$$
Clearly, $X_\cA$ is the closure of the orbit of the point
$(1, \ldots, 1)$. Of course, $\IT^n$ acts trivially on $\IA^n$, by
$$
(s, x) = (s_1, \ldots, s_n, x_1, \ldots, x_n) \mapsto
(s_1 x_1, \ldots, s_n x_n)\ ,
$$
and the embedding $\io:\ X_\cA \hookrightarrow \IA^n$ is 
equivariant with respect to the homomorphism of tori
$$
\vp:\ \IT^d \to \IT^n, \quad t \mapsto (t^{a_1}, \ldots, t^{a_n})\ .
$$

One can then see that $X_\cA$ contains an algebraic
torus $T$ (of dimension $= \rk \cA$) as a dense open subset, and
that there is an action $T \times X_\cA \to X_\cA$ extending the
natural action of $T$ on itself. (See \cite[Ch. I]{KKMS},
\cite[Ch. 5]{GKK} for the converse statement: An affine variety
satisfying the latter condition is of the form $X_\cA$.)
In fact, if $X \subset \IA^n$ is an affine variety defined by
a toric ideal, then $X \cap \IT^n$ is a subgroup of
$\IT^n$. Let $d = \dim X$ and let $\cA = \{a_1, \ldots, a_n\}
\subset \IZ^d$, where the $a_i$ are the columns of a matrix
whose rows form a minimal set of generators of the sublattice
of $\IZ^n$ orthogonal to all $\g \in \IZ^n$ such that
the binomial $x^{\g^+} - x^{\g^-}$ belongs to the ideal of $X$.
Then $\vp: \IT^d \to \IT^n$ (as above) induces an isomorphism onto
$X \cap \IT^n$.

We will sometimes write a  binomial equation $x^{\g^+} - x^{\g^-} = 0$
as $x^{\g} = 1$ (in the ring of Laurent polynomials $\uk[x^\pm]$); of
course, if $I \subset \uk[x]$ is a toric ideal, then $I \cdot \uk[x^\pm]$
defines the torus $T_X$.

We present the following as an exercise using the lemmas
above.

\begin{example}[{\rm Nonexistence of maximal contact in positive
characteristic \cite{Nar}}]
Let $\uk$ denote an algebraically closed field of characteristic $2$.
Let $X \subset \IA^4$ denote the hypersurface
$$
y^2 + x_1^3x_2 + x_2^3x_3 +x_1x_3^7\ =\ 0\ .
$$
Then: (1) The locus of points of $X$ of order $2$ (see Section 3)
is the curve $C$
given by
$$
x_1 = t^{15}\ ,\ \ x_2 = t^{19},\ \ \ x_3 = t^7\ ,\ \ y = t^{32}\ ,
$$
but (2) $C$ lies in no smooth hypersurface. (Hints: (1) The locus
of points of order $2$ is defined by binomial equations. (2) $C$
satisfies no binomial equation of order $1$.)
\end{example}

\subsection{Embedded toric varieties}
Let $N \cong \IZ^n$ denote a lattice, and let $N_\IR$ denote the 
real vector space spanned by $N$. 

A {\it fan} $\Sig$ in $N$ is defined as a (finite) set of {\it cones}
(rational convex polyhedral cones), such that each face of
a cone in $\Sig$ is also a cone in $\Sig$, and the intersection of
two cones in $\Sig$ is a face of each. (See \cite{Ful}, \cite{Oda}.)
Every normal toric variety is associated to a fan \cite[Ch. I]{KKMS}.
A normal toric variety determines the fan; i.e., an
equivariant isomorphism of normal toric varieties corresponds to
an isomorphism of lattices taking one fan to the other 
\cite[Thm. 4.1]{Oda1}. A cone $\s$ in $N_\IR$ corresponds to a
normal affine toric variety $U_\s$. The toric variety associated to
a fan $\Sig$ in $N$ is obtained by glueing together the affine
toric varieties $U_\s$, $\s \in \Sig$; if $\s, \tau \in \Sig$,
then $U_\s \cap U_\tau = U_{\s \cap \tau}$. Clearly, $X = \bigcup U_\s$,
where the union is over the {\it maximal} cones $\s \in \Sig$.
($\s$ is maximal if it is not properly contained in another cone
of $\Sig$.) We will need these general
considerations only in the case of smooth toric varieties.

Let $\s$ denote a cone in $N_\IR$. We define the {\it vertices} of
$\s$ as the (unique) generators $e_\rho$ of the semigroups 
$\rho \cap N$,
for every {\it edge} ($1$-dimensional face) $\rho$ of $\s$. (Then the
vertices of $\s$ generate $\s$.) We will say that $\s$ is {\it regular}
if its vertices $e_{p+1}, \ldots, e_n$ 
can be completed to a basis $e_1, \ldots, e_{p}, e_{p+1},
\ldots, e_n$ of the lattice $N$.

A regular cone $\s$ of dimension $n-p$ in $N_\IR$ corresponds to a 
smooth affine toric variety $U_\s \cong \IT^p \times \IA^{n-p}$;
$U_\s$ can be realized as the subvariety of $\IA^{2p} \times \IA^{n-p}$
(with affine coordinates $(u,x) = (u_1,\dots, u_p, x_1,\dots, x_p,
x_{p+1},\ldots, x_n)$) defined by the binomial equations
$$
u_i x_i\ =\ 1\ , \quad i\ =\ 1, \dots, p\ .
$$

A fan $\Sig$ is {\it regular} if every cone $\s \in \Sig$ is regular.
Smooth toric varieties correspond to regular fans.
Consider a lattice $N \cong \IZ^n$. Let $\Sig \subset N_\IR$ denote a
regular fan, and let $M = M(\Sig)$ denote the smooth toric variety over
$\uk$ determined by $\Sig$. Consider cones $\s, \tau \in \Sig$, of
dimensions $n-p,\ n-q$, respectively. Choose bases $e_1,\ldots, e_n$
and $f_1,\ldots, f_n$ of $N$ such that $e_{p+1}, \ldots, e_n$ and
$f_{q+1}, \ldots, f_n$ are the vertices of $\s$ and $\tau$ 
(respectively).
We can write
$$
e_i\ =\ \sum_{j=1}^n a_{ij} f_j\ , \quad i\ =\ 1, \dots, n\ ,
$$
where $(a_{ij})$ is an integer matrix with determinant $\pm 1$.

Corresponding to the bases $\{e_i\}$ and $\{f_j\}$, we can realize
$U_\s \hookrightarrow \IA^{2p} \times \IA^{n-p}$ with coordinates
$(u,x)$ as above, and $U_\tau \hookrightarrow \IA^{2q} \times \IA^{n-q}$
with analogous coordinates $(v,y) = (v_1,\dots, v_q, y_1,\dots, y_n)$,
in such a way that the transformation from $x$- to $y$-coordinates in
the overlap $U_\s \cap U_\tau$ is given by the Laurent monomial mapping
$$
y_j\ =\ x^{a_j}\ ,\quad \mbox{ where }\ a_j\ =
\ (a_{1j}, \ldots, a_{nj})\ .
$$

A binomial equation $x^{\g} = 1$, $\g \in \IZ^n$, in the 
$x$-coordinates of $U_\s$ becomes $y^{\de} = 1$, where
$\g_i = \sum_{j=1}^n a_{ij} \de_j$, $i = 1, \dots, n$, in $U_\tau$. 
In other words, if $\la:\ N \to \IZ$ is the $\IZ$-linear function
determined by $\la(e_i) = \g_i$, $i = 1, \dots, n$, then $x^{\g} = 1$
in $U_\s$ becomes $y^{\de} = 1$ in $U_\tau$, where $\de_j = 
\la(f_j)$, $j = 1, \dots, n$. This means:

\begin{lemma}
A toric hypersurface
(not necessarily normal) in $M = M(\Sig)$ corresponds to (the
restriction to $\Sig$ of) a linear function $\la:\ N \to \IZ$.
\end{lemma}

In general, if $X$ is a toric subvariety of $M$, then there are
finitely many linear functions $\la_k:\ N \to \IZ$, such that
each $U_\s$ admits a closed embedding in $\IA^{2p} \times \IA^{n-p}$,
for some $p$, with coordinates $(u,x)$ as above, in such a way that
$X \cap U_\s$ is defined by the toric ideal generated by 
$u_j x_j = 1$, $j = 1, \ldots, p$, and $x^{\g_k} = 1$, where
$\g_{ki} = \la_k (e_i)$, $i = 1, \dots, n$, for each $k$.

\section{Equimultiple locus of a binomial}
\subsection{Perfect fields and order of vanishing}
Let $\uk$ denote a field. Let $\um$ be a maximal ideal of
$\uk[x] = \uk[x_1, \ldots, x_n]$. Let $f \in \uk[x]$. The {\it
order} $\mu_{\um} (f)$ of $f$ {\it at} $\um$ is defined as the order of
$f$ as an element of the local ring $R := \uk[x]_{\um}$; i.e.,
the largest $l \in \IN$ such that $f \in \um_R^l$, where $\um_R$
denotes the maximal ideal $\um \cdot R$ of $R$. (The maximal ideals
$\um$ of $\uk[x]$ are the (closed) points $a$ of n-dimensional affine
space $\IA^n = \IA^n_{\uk}$ over $\uk$; we will write $\mu_{\um} (f)$
or $\mu_a(f)$, indifferently.)

There is a second natural notion of order at $a$, corresponding to
the definition above, after extension to the {\it residue field}
$\IF := R / \um_R$. For each $i = 1, \ldots, n$, let $x_i(a) \in \IF$
denote the image of $x_i \in \uk[x]$ in $\IF$; and set $x(a) = 
(x_1(a), \ldots, x_n(a)) \in \IF^n$. Let $\mu_{x(a)}(f)$ denote the 
{\it order at} $x(a)$ of $f$ as an element of $\IF[x]$; i.e., the
largest $l \in \IN$ such that $f \in (x-x(a))^l$, where $(x-x(a))$ 
denotes the maximal ideal of $\IF[x]$ generated by $x_i - x_i(a)$,
$i = 1, \ldots, n$.

We will show that the two notions of order above coincide for
all maximal ideals $\um$ of $\uk[x]$ if and only if $\uk$ is
perfect. 

Let $\um$ be a maximal ideal of $\uk[x] = \uk[x_1, \ldots, x_n]$.
We use the notation above. If $f(x) \in \uk[x]$, let $\of (x) \in
\IF[x]$ denote the image of $f(x)$ by the injective homomorphism
$\uk[x] \hookrightarrow \IF[x]$; then $f(x) \in \um$ if and only
if $\of (x) \in (x-x(a))$. Therefore, there is an induced homorphism
$\tau_{\um} :\ R \to \IF[x]_{(x-x(a))}$. 

The completion of $\IF[x]_{(x-x(a))}$ can be identified with the ring of
formal power series $\IF\lbr X \rbr = \IF\lbr X_1, \ldots, X_n \rbr$,
where the canonical injection $\IF[x]_{(x-x(a))} \hookrightarrow
\IF\lbr X \rbr$ is induced by $g(x) \mapsto g(x(a) + X)$, $g(x)
\in \IF[x]$. Let 
$T_{\um} : \wR \to \IF\lbr X \rbr$ denote the homomorphism of
completions induced by $\tau_{\um}$. (Then $\mu_{x(a)}(f)
= \mu_{(X)}(T_{\um}f) \geq \mu_a(f)$.)
We will need the following
theorem only beginning with Section 7.

\begin{theorem}
The homomorphism $T_{\um}:\ \wR \to \IF\lbr X \rbr$ is an
isomorphism, for all maximal ideals $\um$ of $\uk[x]$, if and
only if $\uk$ is perfect.
\end{theorem}

\begin{proof}
Let $\um$ denote a maximal ideal of $\uk[x]$.
For all $k \in \IN$, $\tau_{\um}$ induces a homomorphism of
finite-dimensional $\IF$-algebras,
$$
T_{\um}^k :\ \frac{\um_R^k}{\um_R^{k+1}} \to \frac{(X)^k}
{(X)^{k+1}}
$$
(where $(X) = (X_1, \ldots, X_n)$ is the maximal ideal of 
$\IF\lbr X \rbr$);
$T_{\um}$ is an isomorphism if and only if $T_{\um}^k$ is an
isomorphism for all $k \in \IN$. 

Suppose that $\uk$ is perfect. For each $i = 1, \ldots, n$, let
$P_i(t) \in \uk[t]$ denote the minimal polynomial of $x_i(a)$. Then
$\overline{P_i}(t) = (t - x_i(a))Q_i(t)$, where $Q_i(t) \in \IF[t]$,
and $Q_i(0) = P_i'(x_i(a)) \neq 0$ \cite[Ch.II, \S5]{ZS}. It is not
difficult to see that $(P_1(x_1), \ldots, P_n(x_n))\cdot \uk[x]$
(the ideal in $\uk[x]$ generated by $P_1(x_1), \ldots, P_n(x_n)$)
is radical, and, as a consequence, that 
$\um_R = (P_1(x_1), \ldots, P_n(x_n))\cdot R$.

If $\al \in \IN^n$, let $Q_\al(x) := \prod_{i=1}^n P_i(x_i)^{\al_i}$;
Then
$$
\overline{Q_\al}(x)\ =\ \xi_\al \cdot (x-x(a))^\al\ \mod 
\ ((x-x(a)))^{|\al|+1}\ ,
$$
where $\xi_\al = \prod P_i'(x_i(a))^{\al_i} \neq 0$, and it follows
that each $T_{\um}^k$ is bijective.

Conversely, suppose that $\uk$ is not perfect. Let $\car \uk 
= p \neq 0$. Take $\la \in \uk \backslash \uk^p$ and set
$f(x) = x^p - \la$ (one variable). Then $f(x)$ is irreducible. Let
$\um$ denote the maximal ideal generated by $f(x)$. We use the
notation above. Let $b$ denote the image of $x$ in $\IF$; then
$b^p = \la$, and 
$$
\tau_{\um} (f)\ =\ (b + (x-b))^p - \la\ =\ b^p + (x-b)^p - \la\ =
\ (x-b)^p\ .
$$
Therefore, $\tau_{\um}(\um_R) \subset (x-b)^p$; in particular
$\um_R / \um_R^2 \to (x-b) / (x-b)^2$ is not an isomorphism.
\end{proof}

\begin{remark} The converse direction in the proof of the theorem
shows that the two notions of order introduced above do not coincide
for all maximal ideals $\um$ when $\uk$ is not perfect. (The other
direction in) the statement of the theorem shows that the two notions
coincide for a perfect field.
\end{remark}

\subsection{Equimultiple locus of a binomial}
Let $f(x) \in \uk[x]$, $x = (x_1, \ldots, x_n)$, be a binomial.
Relabelling the variables, let us write
$$
f(x)\ =\ f(u,v)\ =\ u^\al -v^\be\ ,
$$
where $x = (u,v) = (u_1, \ldots, u_k, v_1, \ldots, v_{n-k})$,
$\al \in \IN^k$, $\be \in \IN^{n-k}$, $0 < |\al| \leq |\be|$ and
$\al_i > 0$, $i = 1, \ldots, k$.

Let $S_f(0)$ denote the {\it equimultiple locus} of the origin
for $f$; i.e.,
$$
S_f(0)\ :=\ \{ b \in \IA_{\uk}^n:\ \mu_b(f) = \mu_0(f) \}\ .
$$
(Of course, $\mu_0(f) = |\al|$.) Clearly, if $\car \uk = 0$
and $|\al| \geq 2$, then 
\begin{equation} \label{equieqn}
S_f(0)\ =\ \{ (u,v):\ u = 0,\ \mu_v(v^\be) \geq |\al| \}\ .
\end{equation}

This is not necessarily true in positive characteristic.

\begin{example}
Let $\car \uk = p > 0$. Let $f(u,v) = u^p - v^{p\be}$ (where
$u$ is a single variable). Then $f(u,v) = (u - v^\be)^p$,
so that $\mu_{(\xi, \eta)}(f) = p$ for any $(\xi, \eta)$ such that
$\xi = \eta^\be$.
\end{example}

We will show that (\ref{equieqn}) holds if 
$f$ is not a $p$'th power in
$\uk[u,v]$, where $p = \car \uk > 0$. We will, in fact, need a
slightly more general result (Theorem \ref{equi} following). If $a \in 
\IA_{\uk}^n$, we let $\IF_a$ denote the residue field of $a$, and
write $x(a) = (x_1(a), \ldots, x_n(a))$, where $x_i(a)$ denotes the
image of $x_i$ in $\IF_a$, for each $i$.

\begin{theorem} \label{equi}
Let $p = \car \uk > 0$. Let
$$
f(u,v,w)\ =\ w^\g u^\al - v^\be\ \in\ \uk[u,v,w]\ ,
$$
where $u = (u_1, \ldots, u_k)$, $v = (v_1, \ldots, v_l)$, $w =
(w_1, \ldots, w_{n-k-l})$. 
Let $a \in \IA_{\uk}^n$. 
Assume that $u(a) = 0$,
$\al_i > 0$ ($i = 1, \ldots, k$), and $w_j(a) \neq 0$ 
($j = 1, \ldots, n-k-l$).
Let $d := \mu_a(f) = |\al|$, and let 
$S_f(a) \subset \IA_{\uk}^n\backslash \{w=0\}$ 
denote the equimultiple locus of $a$; i.e.,
$$
S_f(a)\ =\ \{b \in \IA_{\uk}^n\backslash \{w=0\}:\ \mu_b(f) = d \}\ .
$$
If $f$ is not a $p$'th power in $\uk[u,v,w]$ and $d \geq 2$, 
then $S_f(a) \subset \{u = 0\}$. (If fact, $S_f(a) \not\subset
\{u = 0\}$ if and only if $f$ is  a $d$'th power and $d = p^s$.)
\end{theorem}

\begin{corollary}
Under the hypotheses of Theorem \ref{equi},
$$
S_f(a)\ =\ \{(u,v,w) \in \IA_{\uk}^n\backslash \{w=0\}:\ u = 0, 
\ \mu_v(v^\be) \geq d \}\ .
$$
\end{corollary}

\begin{lemma} \label{equilemma}
Let $d = qp^t$, where $q \not\equiv 0 \mod p$. If $b \in \IA_{\uk}^n$,
$w_j(b) \neq 0$, for all $j$, and $c \in \IF_b$, then
$$
\mu_{w = w(b) \atop z = 0} \left( w^\g (z + c)^d - w^\g c^d \right)\ =\ p^t
$$
(where $z$ is a single variable).
\end{lemma}

\begin{proof}
\begin{align*}
w^\g (z + c)^d - w^\g c^d\ &=\ w^\g \left(z^{p^t} + c^{p^t}\right)^q - w^\g c^d\\
     &=\ w^\g \left( z^d + q z^{(q-1)p^t} c^{p^t} + 
         \cdots + q z^{p^t} c^{(q-1)p^t} \right)\ .
\end{align*}
The result follows.
\end{proof}

\begin{proof}[Proof of Theorem \ref{equi}]
Consider $b \in S_f(a)$; say $(u(b),v(b),w(b)) = (\xi, \eta, \zeta)$.
Then $\mu_{(\xi, \eta, \zeta)}(f) \geq \mu_b(f) = d$. 
Assume that $\xi \neq 0$. We will get a contradiction.

First suppose that $k > 1$. We can assume that $\xi_1 \neq 0$, where
$\xi = (\xi_1, \ldots, \xi_k)$. Write
\begin{equation*}
f(u,v,w)\ =\ w^\g \left( (u - \xi + \xi)^\al - \xi^\al \right)
             - \left( (v - \eta +\eta)^\be - \eta^\be \right)
             + (w^\g \xi^\al -\eta^\be)\ ,
\end{equation*}
and
\begin{equation*}
\left( (u - \xi + \xi)^\al - \xi^\al \right)\  
         =\ (u - \xi)^\al + \sum_{1 \leq |\de| < |\al|}
             {\al \choose \de}(u - \xi)^{\al - \de}\xi^\de\ . 
\end{equation*}
Consider $\de = (\al_1, 0, \ldots, 0)$. Then
\begin{equation*}
{\al \choose \de}(u - \xi)^{\al - \de} \xi^\de\ 
         =\ \xi_1^{\al_1} (u_2 -\xi_2)^{\al_2} \cdots (u_k -\xi_k)^{\al_k}\ ,
\end{equation*}
so that
\begin{equation*}
\mu_{(\xi, \eta, \zeta)}(f)\ \leq\ \al_2 + \cdots + \al_k\ =
\ d - \al_1\ <\ d\ ;
\end{equation*}
a contradiction.

It remain to consider the case that $k = 1$; i.e., $f = w^\g u^d - v^\be$,
where $u$ is a single variable. Write $d = d_1 p^s$, where $d_1 \not\equiv
0 \mod p$. As above, write
\begin{equation*}
f(u,v,w)\ =\ w^\g \left( (u - \xi + \xi)^d - \xi^d \right)
             - \left( (v - \eta +\eta)^\be - \eta^\be \right)
             + (w^\g \xi^d -\eta^\be)\ .
\end{equation*}
By Lemma \ref{equilemma}, 
the first summand of the right-hand side has order $p^s$
in $u - \xi$; therefore, $\mu_b(f) \leq p^s$. So, if $d_1 \neq 1$,
then $\mu_{(\xi, \eta, \zeta)}(f) < d$; a contradiction.

On the other hand, suppose that $d_1 = 1$; i.e., $d = p^s$. Then
\begin{equation*}
f(u,v,w)\ =\ w^\g (u - \xi)^{p^s}
             - \left( (v - \eta +\eta)^\be - \eta^\be \right)
             + (w^\g \xi^{p^s} -\eta^\be)\ .
\end{equation*}
We can assume that $b \in \{f=0\}$; 
i.e., $\zeta^\g \xi^{p^s} - \eta^\be = 0$. We consider two cases:

(1) $\be_j \not\equiv 0 \mod p^s$, for some $j$. Say that 
$\be_1 \not\equiv 0 \mod p^s$; i.e., $\be_1 = qp^t$, where
$q \not\equiv 0 \mod p$ and $0 \leq t < s$. Then
$$
(v - \eta +\eta)^\be - \eta^\be\ =\ (v_1 - \eta_1 +\eta_1)^{\be_1}
\prod_{i \neq 1} (v_i - \eta_i +\eta_i)^{\be_i} - \eta^\be\ .
$$
By Lemma \ref{equilemma}, 
the coefficient of $(v_1 - \eta_1)^{p^t}$ in this expansion
is nonzero. Therefore, $\mu_{(\xi, \eta, \zeta)}(f) \leq p^t < d$; 
a contradiction.

(2) $\be_j \equiv 0 \mod p^s$, for all $j$. 
Then
\begin{align*}
f(u,v,w)\ &=\ (w - \zeta + \zeta)^\g \left( (u - \xi)^{p^s} + \xi^{p^s}\right)
             - (v - \eta + \eta)^\be\\
         &=\ \left( (w - \zeta + \zeta)^\g - \zeta^\g \right) \xi^{p^s}
             + (w - \zeta + \zeta)^\g (u - \xi)^{p^s}
             - \left( (v - \eta + \eta)^\be - \eta^\be \right) .
\end{align*}
The second summand here has order at least $p^s$ (with respect to 
$(u - \xi, v - \eta, w - \zeta)$), and the third summand also has order at
least $p^s$, since all $\be_j \equiv 0 \mod p^s$. Thus 
$\left( (w - \zeta + \zeta)^\g - \zeta^\g \right)$ has order at least
$p^s$ in $w - \zeta$. If all $\g_j \equiv 0 \mod p^s$, then
$f$ is a $p^s$'th power; a contradiction to the hypothesis of the
theorem. On the other hand, if some $\g_j \not\equiv 0 \mod p^s$,
then we can show that $\left( (w - \zeta + \zeta)^\g - \zeta^\g \right)$ 
has order $< p^s$, using Lemma \ref{equilemma}
as in the case (1); a contradiction.
\end{proof}

\begin{remark}
Let $\uk$ be a field. Let $f(x) \in \uk[x] = \uk[x_1, \ldots, x_n]$.
If $X = (X_1, \ldots, X_n)$, then $f(x+X) \in \uk[x][X]$; say,
$f(x+X) = \sum f_\al(x) X^\al$. For each $\al \in \IN^n$, 
$f_\al(x)$ is called the {\it Hasse derivative} of $f$ of {\it order}
$\al$. Let $f(u,v,w)$ be a binomial as in Theorem \ref{equi}. By 
Corollary 3.5, the equimultiple locus $S_f(a)$ 
has the structure of a sub{\it space} (or subscheme) of 
$\IA_{\uk}^n \backslash \{w=0\}$ given by the
vanishing of all Hasse derivatives of $f$ (or all Hasse derivatives
with respect to 
$(u,v)$) of orders $< d$. (In particular, the order $\mu_a(f)$ is
a Zariski upper-semicontinuous function of $a$.) According to \S3.1
above, the analogous result for general polynomials requires
the hypothesis that $\uk$ be perfect.

\end{remark}

\section{Blowing up and strict transform}
Let $\Sig$ denote a regular fan in a lattice $N \cong \IZ^n$. (See
\S2.2.) Let $\D$ be a cone in $\Sig$, and let $e_\D \in N$ denote the
sum of the vertices of $\D$. We call $e_\D$ the {\it barycentre}
of $\D$. We define the {\it star-subdivision} $\Sig'$ of $\Sig$
as the smallest refinement of $\Sig$ that includes $e_\D$ as
a vertex.

Let $M = M(\Sig)$ denote the smooth
toric variety over $\uk$ corresponding to $\Sig$. 
Let $\la:\ N \to \IZ$ be
a $\IZ$-linear function, and let $X \hookrightarrow M$ denote the
corresponding toric hypersurface in $M$ (Lemma 2.6).

\subsection{Affine case}
$M = \IA^n = \IA^n_{\uk}$. We can assume that $N = \IZ^n$ and that
$\Sig$ is given by the cone generated by the standard basis vectors 
$e_1, \ldots, e_n$. In the affine coordinates $x = (x_1, \ldots, x_n)$
of $\IA^n$, the ideal of $X$ is generated by the binomial $x^\g - 1$
(in $\uk[x^\pm]$),
where $\g = (\g_1, \ldots, \g_n)$ and $\g_i = \la(e_i)$, $i = 1, \ldots,
n$ (Lemma 2.1 and \S2.2).

The $\IT^n$-invariant subspaces of $\IA^n$ are simply the
coordinate subspaces $Z_\D := \{x_i = 0,\ i \in \D\}$, where
$\D \subset \{1, \ldots, n\}$. Fix such $\D$, and consider the
blowing-up $\pi$ of $M' \to M = \IA^n$ with centre $Z_\D$. 
This blowing-up
can be described combinatorially as follows. We identify $\D$
with the face of $\Sig$ spanned by $e_i$, $i \in \D$. Then 
$M' = M(\Sig')$, where $\Sig'$ is the star-subdivision of $\Sig$
determined by $\D$. For each $i \in \D$, let $\s_i \in \Sig'$
denote the cone with vertices $e_j$, $j \in \D \backslash \{i\}$,
and $e_\D$. (The $\s_i$, $i \in \D$, are the maximal cones of $\Sig'$.)
For each $i \in \D$, $U_{\s_i}$ has affine coordinates 
$x = (x_1, \ldots, x_n)$ with respect to which $\pi: U_{\s_i} \to
\IA^n$ is given by the following substitution rules: For each
$j \in \D \backslash \{i\}$, substitute $x_i x_j$ for $x_j$, and
leave the remaining variables $x_j$ unchanged. (For economy of
notation, we use the same symbols for the variables before and
after blowing up.)

Let $\g_\D = \la(e_\D) = \sum_{j \in \D} \g_j$. It is easy to
see that the strict transform $X'$ of $X$ by $\pi$ is defined
in each affine chart $U_{\s_i}$, $i \in \D$, by the toric ideal
generated by $x^{\g'} -1$, where $\g' = (\g'_1, \ldots, \g'_n)$
and $\g'_j = \g_j$, $j \neq i$, and $\g'_i = \g_\D$. In other
words, $X'$ is defined by (the restriction to $\Sig'$ of) the
same linear function $\la:\ N \to \IZ$ as before blowing-up.

\subsection{General case}
Let $\D$ denote a cone in $\Sig$. Then $\D$ determines a smooth
invariant subspace $Z_\D$ of $M = M(\Sig)$: Consider any cone
$\s \in \Sig$ such that $\D$ is a face of $\s$;
choose a basis $e_1, \ldots, e_n$ of $N$ so that
$e_{p+1}, \ldots, e_n$ are the vertices of $\s$. We can realize $U_\s$
as a closed subvariety of $\IA^{2p} \times \IA^{n-p}$ with affine
coordinates $(x,u)$ as in \S2.2. Then $Z_\D \cap U_\s$ is defined
by the equations $u_i x_i = 1$, $i = 1, \ldots, p$, and $x_j = 0$,
$j \in \D$. 

Let $\pi:\ M' \to M = M(\Sig)$ denote the blowing-up with centre
$Z_\D$. The $M' = M(\Sig')$, where $\Sig'$ is the star-subdivision
of $\Sig$ determined by $\D$. The strict transform $X'$ of $X$ by 
$\pi$ is the toric hypersurface of $M'$ determined by (the 
restriction to $\Sig'$ of) the linear function $\la$.

\section{Equivariant resolution of singularities of a toric hypersurface}
In this section, we prove Theorem 1.4; we translate Theorem 1.4 into
a purely combinatorial statement (Theorem 5.3 below) which we prove
using a simple algorithm. Our desingularization algorithm for a toric
hypersurface corresponds to the ``combinatorial
part'' of the desingularization algorithm of \cite{BMinv}. (See
\cite[p. 260]{BMinv}, \cite[\S3.3(1)]{BMda1} and also 
\cite[\S4]{BMihes}, as well as Step 2(b) in the proof of Theorem 8.5
below.) Throughout this section, we will use the language
of Section 4 without further notice.

\subsection{Affine case}
We first consider the special case of an affine toric hypersurface
$X \subset \IA^n = \IA_{\uk}^n$. 
The toric ideal corresponding to 
$X$ is generated by a binomial $x^\g - 1$; i.e., by a binomial
$f(x) = x^{\g^+} - x^{\g^-}$ in the affine coordinates $x = (x_1, 
\ldots, x_n)$ of $\IA^n$, where $\g \in \IZ^n$ (Lemma 2.1). 
We can assume that $|\g^-| \leq |\g^+|$; then $d := |\g^-|$ is
the maximum order $\mu_0(f)$ of $f$ (or $X$). (The maximum order
of a binomial is, of course, taken at the origin.) 

The locus of
maximum order $d$ of $X$ is the equimultiple locus $S_f(0)$
of $0$. By Corollary 3.5, 
$$
S_f(0)\ =\ \{x \in \IA^n:\ x_i = 0\ \mbox{ if }\ \g_i < 0,
\ \mbox{ and }\ \mu_x(x^{\g^+}) \geq d\}\ .
$$
Let $Z_\D := \{x_i = 0,\ i \in \D\}$, where
$\D \subset \{1, \ldots, n\}$. Then $Z_\D \subset
S_f(0)$ if and only if
\begin{equation}
\begin{array}{cc}
i \in \D \ \mbox{ if }\ \g_i < 0 
&(\mbox{i.e., }\ \g^-_\D = d\, )\ ,\hfill\\
\g^+_\D \geq d = \g^-_\D 
&(\mbox{i.e., }\ \g_\D \geq 0\, )\hfill
\end{array}
\end{equation}
(recall that $\g_\D := \sum_{i \in \D} \g_i$; likewise for
$\g^+, \g^-$), and 
\begin{equation}
S_f(0)\ =\ \bigcup Z_\D\ ,
\end{equation}
where the union is over the {\it minimal} subsets $\D$ of
$\{1, \ldots, n\}$ satisfying (5.1).
In other words, the union in (5.2) is over the
subsets $\D$ of $\{1, \ldots, n\}$ such that
\begin{equation}
\begin{array}{cc}
\g^-_\D = d\ ,&\\
0 \leq \g_\D < \g_i\ ,& \mbox{for all }\ i \in \D\ 
\mbox{ such that }\ \g_i \geq 0\ .\hfill
\end{array}
\end{equation}
(In particular, $\g_i \neq 0$, for all $i \in \D$.)

Consider the blowing-up of $\IA^n$ 
with centre $Z_\D \subset S_f(0)$. Let
$X'$ denote the strict transform of $X$. Then (in the notation of
\S4.1) $X'$ has order $\leq d$ at the origin of each chart 
chart $U_{\s_i}$, $i \in \D$, and therefore at every point. 

Now suppose that $Z_\D$ is a component of $S_f(0)$. It follows
from (5.3) that
$X'$ has order $< d$ throughout each $U_{\s_i}$ where $\g_i < 0$,
and, in a chart $U_{\s_i}$, $\g_i > 0$, if $\mu_{f'}(0) = d$, where
$f'(x) = x^{(\g')^+} - x^{(\g')^-}$, then 
$$
d\ \leq\ |(\g')^+|\ <\ |\g^+|\ .
$$

It follows that we can reduce the order over every chart by a finite
number of analogous blowings-up. At each step, the centre of blowing-up
extends to a global closed smooth invariant subspace, as we will
observe in \S5.2 below. The main point is the following: For any
$\D \subset \{1, \ldots, n\}$ (i.e., for any face $\D$ of $\Sig$,
in the language of \S4.1), set
\begin{align*}
d_\D &:= \min\{\g^-_\D, \g^+_\D\}\ ,\\
\Om_\D &:= \max\{\g^-_\D, \g^+_\D\}\ .
\end{align*}
(In particular, $d_{\Sig} = d$. Then $Z_\D$ is a component of $S_f(0)$
above if and only if
\begin{align*}
d_\D &= d\ ,\\
d_{\D_1} &< d\ ,\mbox{ for every proper subface } \D_1 \mbox{ of } \D\ .
\end{align*}
Moreover, if $\Sig'$ is the star-subdivision of $\Sig$ determined
by such a face $\D$, then
$$
(d_\s, \Om_\s)\ <\ (d_{\Sig}, \Om_{\Sig})
$$ 
(with respect to the lexicographic ordering of pairs), for every
cone $\s \in \Sig'$.

\subsection{General case}
We now consider a general toric hypersurface $X \subset M = M(\Sig)$
over $\uk$, as in Section 4. For any cone $\s 
\in \Sig$, write
$$
\g_\s^+ := \sum_{e \in \s \atop \la(e) > 0} \la(e)\ , \qquad
\g_\s^- := \sum_{e \in \s \atop \la(e) < 0} (-\la(e))
$$
(where ``$e \in \s$'' means that $e$ is a vertex of $\s$). Set
\begin{align*}
d_\s &:= \min\{ \g_\s^-, \g_\s^+ \}\ ,\quad \s \in \Sig\ ,\\
d(\la, \Sig) &:= \max_{\s \in \Sig} d_\s\ .
\end{align*}

\begin{definitions}
We will say that a cone $\D \in \Sig$ (or the star-subdivision of
$\Sig$ determined by $\D$) is
\begin{enumerate}
\item {\it admissible}\, if $d_\D = d(\la, \Sig)$.
\item {\it minimal}\, if it is admissible and $d_{\D_1} < d(\la, \Sig)$,
for every proper subface $\D_1$ of $\D$.
\end{enumerate}
\end{definitions}

Clearly, $d(\la, \Sig)$ is the maximum order of $X$ (if
$d(\la, \Sig) > 0$).
We recall that the centre of blowing up $Z_\D$ corresponding
to the star-subdivision of $\Sig$ determined by a cone $\D$ is
a smooth closed invariant subspace of $M$. (Therefore, $Z_\D$
simultaneously has only normal crossings with respect to all
smooth invariant subspaces of $M$.) If $\D$ is admissible, then
$X$ assumes its maximum order $d(\la, \Sig)$ at each point of $Z_\D$
(so $X$ is normally flat along $Z_\D$). If $\D$
is minimal, then $Z_\D$ is also a component of the locus of points
of $X$ of maximum order. In particular, we obtain the following
version of Theorem 1.5.

\begin{proposition}
The locus of maximum order (the {\it equimultiple locus} $S$) of
a toric hypersurface $X \subset M$
has only normal crossings. Moreover, $S$ is a union of global
smooth closed $T_M$-invariant components.
\end{proposition}

The following theorem is a combinatorial restatement of Theorem 1.4.

\begin{theorem}
There is a finite succession of admissible star-subdivisions,
\begin{equation*}
\Sig\ =\ \Sig_0\ , \Sig_1\ , \ldots,\ \Sig_t\ ,
\end{equation*}
such that $d(\la, \Sig_t) = 0$.
\end{theorem}

\begin{remark}
For each $l \geq 0$, let $X_l$ denote the toric hypersurface 
in $M_l := M(\Sig_l)$ determined by the function $\la$
on $\Sig_l$. Then $X_{l+1}$ is the strict transform of $X_l$
by the blowing-up $M_{l+1} \to M_l$ determined by the star-subdivision
$\Sig_{l+1}$ of $\Sig_l$. The condition $d(\la, \Sig_t) = 0$ means
(in the notation of \S\S2.2, 4.2) that, in each chart $U_\s$ of $M_t$,
$X_t$ is defined by a binomial $x^\g -1$, where $\g_i \geq 0$ for
all $i$. Therefore, $X_t$ is smooth and simultaneously has only
normal crossings with respect to the collection of exceptional
hypersurfaces (all given by coordinate subspaces in the local charts).
\end{remark}

Our remarks in \S5.1 show that
we can simply use the following combinatorial
algorithm to prove Theorem 5.3.

\begin{algorithm}
For each $l \geq 0$, let $\Sig_{l+1}$ be any minimal star-subdivision 
of $\Sig_l$.
\end{algorithm}

\begin{questions}
Let us say that a sequence of star-subdivisions, $\Sig = \Sig_0,
\Sig_1, \ldots, \Sig_t$, is {\it resolving} if $d(\la, \Sig_t) = 0$.

(1) Does any resolving sequence of star-subdivisions have length
greater than or equal to the length of some resolving sequence of
minimal star-subdivisions?

(2) It is not difficult to give an example of an affine toric
hypersurface for which there are resolving sequences of minimal
star-subdivisions of different lengths. Do we get a resolving
sequence of shortest length by taking, at each step, 
the star-subdivision based on any minimal cone $\D$ with the
smallest value of $\Om_\D = \max\{ \g_\D^-, \g_\D^+ \}$?
\end{questions}

\section{Standard basis of a toric or binomial ideal}
Our goal in this section is to show that a standard basis of 
a toric or more general binomial ideal is given by binomials. 
(Compare \cite{Stu1}.) Let $\uk$ denote a field and let
$\uk[x] = \uk[x_1,\ldots, x_n]$.

\subsection{Binomial ideal} Let $I$ denote an ideal in $\uk[x]$.

\begin{definition} 
We say that $I$ is a {\it binomial ideal} if:
\begin{enumerate}
\item $I$ is generated by binomials $x^{\g^+} - x^{\g^-}$, $\g \in \IN$.
\item $\uk[x]/I$ contains no nilpotents.
\item If $f(x) \in \uk[x]$ and $x_i f(x) \in I$ (for some $i=1,\ldots,n$),
then $f(x) \in I$.
\end{enumerate}
\end{definition}

A toric ideal is a special case of a binomial ideal. (If $X$ is the
affine toric variety $V(I)$ corresponding to a toric ideal $I$, then
property (3) above expresses the density of $T_X = X \cap \IT^n$ in
$X$.) Note that we use ``binomial ideal'' in a more restrictive sense
than \cite{ES}.

If $I \subset \uk[x]$ is a binomial ideal, we call $X = V(I) \subset
\IA^n$ an {\it affine binomial variety}. Affine binomial varieties
share many of the properties of affine toric varieties.

\subsection{Distinguished point} An affine binomial variety $X \subset
\IA^n = \IA_{\uk}^n$ has a {\it distinguished point} $a$: Let $I \subset
\uk[x]$ denote the binomial ideal corresponding
to $X$. (After permuting the variables if necessary) we can assume
that $x_{n-m+1}, \ldots, x_n$ are precisely the variables $x_i$ that
vanish nowhere on $X$. Let us relabel the variables $(x_1, \ldots, x_n)$
as $(x,y) = (x_1, \ldots, x_{n-m}, y_1, \ldots, y_m)$. Then it is
easy to see that $I$ (as an ideal in $\uk[x,y^\pm]$) 
is generated by binomials of the form
\begin{equation*}
\begin{array}{cc}
1 - y^\g\ ,& \ \g \in \IZ^m\ ,\hfill\\
x^\al - x^\be y^\g\ ,& \ \al, \be \in \IN^{n-m},\ \g \in \IZ^m\ ,
\ 1 \leq |\al| \leq |\be|\ .\hfill
\end{array}
\end{equation*}
Let $a$ denote the point $(\uzero, \uone)$, where 
$\uzero = (0, \ldots, 0)
\in \uk^{n-m}$, $\uone = (1, \ldots, 1) \in \uk^m$.

In the toric case, 
the distinguished point $a$ belongs to the unique closed orbit of
$X$ (in particular, $a$ belongs to the closure of every orbit).
The closed orbit of $X$ is an algebraic torus isomorphic to $T_Y$,
where $Y \subset \IA^m$ is the toric variety defined by the ideal
generated by all $1 - y^\g \in I$, and $a$ is, in fact, its identity
element \cite[Ch. 3]{Ful}

Let $X \subset \IA^n$ be an affine binomial variety, as above.
Consider any affine subspace (i.e., coordinate subspace)
$\IA^q$ of $\IA^n$ ($q \leq n$), and
let $\IT^q \subset \IA^q$ denote the standard torus. Then (as in the
toric case), $X \cap \IT^q$ is a subgroup of $\IT^q$. In particular,
$X \cap \IT^q$ is smooth. 

It follows that (in the notation above), the affine variety $Y \subset
\IA^m$ defined by the ideal generated by all $1 - y^\g \in I$ is a
smooth affine binomial variety (toric, if $I$ is a toric ideal).

\begin{lemma}
If $X \subset \IA^n$ is an affine binomial variety, then $X \cap \IT^n$
is dense in $X$.
\end{lemma}

\begin{proof}
Let $I \subset \uk[x_1,\ldots, x_n]$ be the binomial ideal of $X$.
If $f(x) \in \uk[x]$ vanishes on $X \cap \IT^n$, then 
$x_1\cdots x_n f(x)$ vanishes on $X$, so that
$x_1^r\cdots x_n^r f(x)^r \in I$, for some positive integer $r$. 
Therefore,
$f(x)^r \in I$, by Definition 6.1(3), so that $f(x)\in I$, by (2).
\end{proof}

\begin{lemma}
Let $X \subset \IA^n$ be an affine binomial variety. Then the 
distinguished point $a \in X$ belongs to the closure of $X \cap \IT^q$,
for every affine subspace $\IA^q \subset \IA^n$ such that $X \cap \IT^q
\neq \emptyset$.
\end{lemma}

\begin{proof}
Let $\IA^q$ be an affine subspace of $\IA^n$, and let $I \subset
\uk[x,y^\pm]$ denote the ideal of $X$ (using the notation above).
(After permuting the variables $x_1, \ldots x_{n-m}$ if necessary)
we can assume that $x = (u_1 \dots, u_{n-q}, v_1, \ldots, v_{q-m})$
($q \geq m$) and that $\IA^q = \{u=0\}$. (Otherwise, $X \cap \IT^q
= \emptyset$.) Let $W \subset \IA^n$ denote the open subset defined
by $v_i \neq 0$ and $y_j \neq 0$, for all $i = 1, \ldots, q-m$ and
$j = 1, \ldots, m$. Then $X \cap \IT^q = X \cap \{u=0\} \cap W$ is
a smooth subvariety of $\{u=0\} \cap W$ defined by those binomials
in $I$ which involve only the variables $v_i, y_j$.

Let $Z \supset X$ denote the binomial variety defined by the latter
binomials. (In particular, $X \cap \{u=0\} \cap W 
= Z \cap \{u=0\} \cap W$.) $Z \cap W$ is a smooth subvariety of $W$.
Since $X \cap \IT^n$ is dense in $X$, by Lemma 6.2, $a \in
\overline{X \cap \IT^n} \subset \overline{Z \cap W}$; therefore,
$a$ is in the closure of $Z \cap \{u=0\} \cap W$, as required.
\end{proof}

\subsection{Diagram of initial exponents} Let $A$ be a commutative 
ring with identity. Let $A\lbr x \rbr = A\lbr x_1, \ldots, x_q \rbr$.
If $\al = (\al_1, \ldots, \al_q) \in \IN^q$, put $|\al| = \al_1
+ \cdots + \al_q$. We will use the total order on $\IN^q$ that is
given by the lexicographic ordering of $(q+1)$-tuples 
$(|\al|, \al_1, \ldots, \al_q)$. Let $F = \sum_{\al \in \IN^q}
F_\al X^\al \in A\lbr x \rbr$, where $x^\al := x^{\al_1} \cdots 
x^{\al_q}$. Let $\supp F := \{\al: F_\al \neq 0 \}$. The {\it initial
exponent} $\exp F$ is the smallest element of $\supp F$. ($\exp F
:= \infty$ if $F = 0$.) If $F \neq 0$ and $\al = \exp F$, then
$F_\al X^\al$ is called the {\it initial monomial} $\mon F$ of $F$.

Let $I$ be an ideal in $A\lbr x \rbr$. The {\it diagram of
initial exponents} $\frN(I) \in \IN^q$ is defined as
$$
\frN(I)\ :=\ \{\exp F:\ F \in I \backslash \{0\} \}\ .
$$
Clearly, $\frN(I) + \IN^q = \frN(I)$. It follows that there
is a smallest finite subset $\frV$ of $\frN(I)$ (the {\it
vertices} of $I$) such that $\frN(I) = \frV + \IN^q$.

\subsection{The vertices of the diagram of initial exponents
of a toric or binomial ideal are represented by binomials} 
In the remainder of Section 6, we assume that $I$ is any binomial
ideal in $\uk[x,y] = \uk[x_1, \ldots, x_{n-m}, y_1, \ldots, y_m]$
and that $X = V(I)$ has nonempty intersection with $\{0\} \times
\IT^m \subset \IA^{n-m} \times \IA^m$. (For desingularization of
toric varieties, we will be interested only the in special case
that $I$ is a toric ideal and $a := (\uzero, \uone)$ is the 
distinguished point, as above.) We write binomials in $I$ in the
form $x^\al - x^\be y^\g$,
where $\g \in \IZ^m$; i.e., we consider $I$ as an ideal in 
$\uk[x,y^{\pm}]$. Then any binomial in $I$ either involves $y$
alone, or is of the form $x^\al - x^\be y^\g$, where $\al, \be$
are nonzero elements of $\IN^{n-m}$.

Let $J \subset \uk[y^\pm]$ denote the ideal generated by all
binomials $1 - y^\g \in I$, and
let $B$ denote the quotient ring $\uk[y^{\pm}]/J$. If $f \in
\uk[x,y^{\pm}]$, let $f_J \in B[x]$ denote the element induced
by $f$. Let $I_J$ denote the ideal in $B[x]$ induced by $I$.

We consider $B[x]$ as a subring of the ring of formal power
series $B\lbr x \rbr$, and write $\frN(I_J) \subset \IN^{n-m}$
for the diagram of
initial exponents of $I_J\cdot B\lbr x \rbr$. 

\begin{lemma} \label{6.4}
{\rm (1)} If $f \in I$ and $\supp f_J \subset \IN^{n-m}\backslash
\frN(I_J)$, then $f \in J\cdot \uk[x,y^{\pm}]$.

{\rm (2)}
If $\al$ is a vertex of $\frN(I_J)$, then $\al = \exp G_J$, where
$G \in I$ is a binomial of the form $x^{\al} - x^{\be}y^{\g}$,
with $\be \in \IN^{n-m}$, $\al < \be$,
and $\g \in \IZ^m$.
\end{lemma}

\begin{proof}
The assertion (1) is obvious. For (2), let
\begin{align*}
\frN\ :=\ \{ \exp f_J \in \IN^{n-m}:\ &f = x^\al - x^\be y^\g \in I,\\
	    &\al, \be \in \IN^{n-m},\ 0 < \al < \be,\ \g \in \IZ^m \}\ ,
\end{align*}
The assertion (2) means that $\frN(I_J) = \frN$.
Clearly, $\frN \subset \frN(I_J)$. Consider $f \in I$. Then
we can write
\begin{equation} \label{eqn6.1}
f\ =\ \sum_{i=1}^t g_i f_i\ \mod\ J\cdot \uk[x,y^{\pm}]\ ,
\end{equation}
where each $g_i \in \uk[x,y^{\pm}]$ and each $f_i$ is a binomial 
\begin{equation*}
f_i\ =\ x^{\al^i} - x^{\be^i} y^{\g^i}\ \in\ I\ ,
\end{equation*}
in which $\al^i, \be^i \in \IN^{n-m}$, 
$\al^i < \be^i$, and $\g^i \in \IZ^m$. Set
\begin{equation*}
\de\ :=\ \min_i \{ \exp g_{iJ} + \al^i \}\ .
\end{equation*}

If $\de = \infty$, then $f_J = 0$; i.e., $f \in J\cdot \uk[x,y^{\pm}]$.
Suppose that $\de < \infty$. Then $\de \leq \exp f_J$.  
If $\exp f_J = \de$, then $\exp f_J \in \frN$. Suppose that
$\de < \exp f_J$. Let
\begin{equation*}
\La\ :=\ \{ i:\ \exp g_{iJ} + \al^i = \de \}\ .
\end{equation*}
If $i \in \La$, then $\mon g_{iJ} = \mu_i x^{\de - \al^i}$,
where $\mu_i$ is a nonzero element of $B$; take $\la_i \in
\uk[y^\pm]$ such that $\mu_i$ is induced by $\la_i$. Then
$\sum_{i \in \La} \la_i \in J$.  We can assume that 
$\La = \{1, \ldots, s\}$, where $2 \leq s \leq t$. Write
\begin{align*}
g'_i &:= g_i - \la_i(y)x^{\de -\al^i}\ ,\quad
i = 1,\ldots,s\ ,\\
g'_i &:= g_i\ , 
\phantom{ - \la_i(y)x^{\de -\al^i}\, ,}\quad
i = s+1, \ldots, t\ .
\end{align*}
Then, modulo $J\cdot\uk[x,y^\pm]$, 
\begin{align}
f &= \sum_{i=1}^t g'_i f_i 
+ \sum_{j=1}^s \la_j(y)x^{\de -\al^j}f_j \notag \\
&= \sum_{i=1}^t g'_i f_i 
+ \sum_{j=2}^s \la_j(y)(x^{\de -\al^j}f_j - x^{\de -\al^1}f_1)\ .
\end{align}
For each $j =2, \ldots, s$, 
\begin{equation*}
x^{\de -\al^j}f_j - x^{\de -\al^1}f_1\ =\ 
\pm y^{\s^j} (x^{\tal^j} - x^{\tbe^j} y^{\tg^j})\ ,
\end{equation*}
where $\s^j, \tg^j \in \IZ^m$ and
\begin{align*}
\tal^j &= \min \{\de-\al^j+\be^j,\ \de-\al^1+\be^1\}\ ,\\
\tbe^j &= \max \{\de-\al^j+\be^j,\ \de-\al^1+\be^1\}\ .
\end{align*}
Note that $\exp g'_{iJ} > \exp g_{iJ}$, $i = 1, \ldots, s$ and
$\tal^j > \de$, $j = 2, \ldots, s$. If 
$\de-\al^j+\be^j = \de-\al^1+\be^1$, for some $j$, then
$x^{\tal^j} - x^{\tbe^j} y^{\tg^j} = x^{\tal^j}(1 - y^{\tg^j}) \in
I$, so that $1 - y^{\tg^j} \in J$ (by Definition 6.1(3)).
Therefore, (6.2) is a new representation of $f$ of the same form
as (6.1) and, if $\de'$ denotes the analogue of $\de$ for this
new representation, then
$$
\de < \de' < \exp f_J\ .
$$
By induction, $f$ has a representation of the form (6.1) with
$\de = \exp f_J$, and the result follows.
\end{proof}

\subsection{The standard basis is given by binomials} 
Let $\al^i$, $i = 1, \ldots, s$, denote the vertices of $\frN(I_J)$.
We associate to $\al^1, \ldots \al^s$ a partition of 
$\IN^{n-m}$: Set $\D_i := (\al^i + \IN^{n-m}) 
\backslash \cup_{j=1}^{i-1} \D_j$, $i = 1, \ldots, s$,
and put $\D_0 := \IN^{n-m} \backslash \cup_{i=1}^s \D_i$.
Cleary, $\frN(I_J) = \cup_{i=1}^{s} \D_i$ and 
$\D_0 = \IN^{n-m} \backslash \frN(I_J)$. Let $A$ denote the
ring of Laurent polynomials $\uk[y^\pm]$.
By Lemma 6.4(2), for each $i = 1, \ldots, s$, there is a binomial 
$G^i = x^{\al_i} - x^{\eta_i} y^{\zeta_i}$
in $I \subset \uk[x,y^\pm] = A[x]$ that represents $\al^i$.

\begin{theorem}[{\rm Hironaka division}] \label{6.5}
{\rm (1)} For each $F \in A[x] \subset A\lbr x \rbr$, there are
unique $Q_i, R \in A[x]$ such that $\al^i + \supp Q_i \subset
\D_i$ ($i = 1, \ldots, s$), $\supp R \subset \D_0$, and 
$$
F\ =\ \sum_{i=1}^s G^i Q_i + R\ .
$$

{\rm (2)} Moreover, if $F = x^\de H(y)$, where $\de \in \IN^{n-m}$
and $H(y)$ is a Laurent monomial, then the remainder $R$ has
the same form.

{\rm (3)} For each $i = 1, \ldots, s$, there is a unique binomial
$$
F^i\ =\ x^{\al^i} - x^{\be^i}y^{\g^i}\ \in\ I,
$$
where $\be^i \in \IN^{n-m} \backslash \frN(I_J)$, $\al^i < \be^i$,
and $\g^i \in \IZ^m$.
\end{theorem}

\begin{remark} \label{6.6}
It follows from Theorem 6.5(1) and Lemma 6.4(1), that $I$ is
generated by the binomials $F^i = x^{\al^i} - x^{\be^i}y^{\g^i}$, 
modulo $J\cdot A[x]$. We call $\{ F^i \}$ the
{\it standard basis} of $I$ modulo $J\cdot A[x]$.
\end{remark}

\begin{proof}[Proof of Theorem \ref{6.5}]
It is enough to prove (1) in the case that $F = x^\de H(y)$, where
$\de \in \IN^{n-m}$ and $H(y) \in \uk[y^\pm]$ is a Laurent monomial.
We will show that, under this assumption, we obtain (1) with
the stronger conclusion of (2).

Consider any $F \in A[x]$. Clearly, there
are unique $Q_i(F),\, R(F) \in A[X]$ such that
$\al^i + \supp Q_i(F) \in \D_i$, $\supp R(F) \in \D_0$, and
$F = \sum_{i=1}^s X^{\al^i} Q_i(F) + R(F)$. Set
\begin{equation} \label{eqn7.3}
E(F)\ :=\ F - \sum G^i Q_i(F) -  R(F)\ =\ \sum (X^{\al^i} 
- G^i) Q_i(F)\ .
\end{equation}
Put $E^0(F) = F$ and $E^j(F) = E(E^{j-1}(F))$, $j \geq 1$.

Now assume that $F = x^\de H(y)$, as in (2), where $\de \neq 0$
and $H \neq 0$. There are two cases:
\begin{enumerate}
\item[(i)] If $\de \not\in \frN(I_J)$, then $Q_i(F) = 0$,
($i = 1, \ldots, s$), $R(F) = F$ and $E(F) = 0$. In this case,
$F = \sum G^i Q_i + R$, where $Q_i = Q_i(F) = 0$, $i = 1, \ldots, s$,
and $R = R(F) = F$.
\item[(ii)] Suppose that $\de \in \frN(I_J)$. 
Then there is a unique $i$ such 
that $\de \in \D_i$; i.e., $\de = \al_i + \ep$, 
where $\ep \in \IN^{n-m}$.
Therefore, 
$$
x^\de H(y)\ =\ x^\ep(x^{\al^i} - x^{\eta^i}y^{\zeta^i})H(y) 
+ x^{\eta^i+\ep}y^{\zeta^i}H(y)\ .
$$
So $R(F) = 0$ and $E(F) = x^{\eta^i+\ep}y^{\zeta^i}H(y)$.
Moreover, 
$$
\exp F = \al^i + \ep < \eta^i + \ep = \exp E(F)\ .
$$
\end{enumerate}

In case (ii), we can continue to divide. We claim there exists
$t \in \IN \backslash \{ 0\}$ such that: 
\begin{enumerate}
\item[(a)] For each $j = 0 \ldots, t$, $E^j(F)$ is a nonzero
element of $A[x]$ of the form $x^{\de_j} H_j(y)$, where
$\de_j \in \IN^{n-m}$ and $H_j(y)$ is a Laurent monomial;
\item[(b)] $E^{t+1}(F) = 0$.
\end{enumerate}

Suppose that this is not so. Then, for all $j = 0,1,\ldots$,
$R(E^j(F)) = 0$ and $E^j(F) = x^{\de_j} H_j(y) \in A[x]$,
where $\de_{j+1} > \de_j$. Therefore, $F = x^\de H(y) \in
I\cdot A\lbr x \rbr$, and it follows that
$x^\de H(y)$ belongs to the ideal generated by $I$ in the localization
of $A[x] = \uk[x,y^\pm]$ at the point $(\uzero,\uone)$ (by faithfull
flatness of completion of a Noetherian local ring).
Then $1$ belongs to this ideal (by Definition 6.1(3)); a
contradiction.

Assertions (1) and (2) follow.

It follows from (2) that, if $F = x^{\al^i}$ (where $1 \leq i \leq s$),
then the remainder $R$ is a nonzero element of $A[x]$ of the form
$x^{\be^i} y^{\g^i}$, where $\be^i \in \IN^{n-m}$, $\be^i > \al^i$,
and $\g^i \in \IZ^{m}$. (3) follows.
\end{proof}

\begin{remark}
Let $V(J) \subset \IA^m$ denote the affine binomial variety determined
by $J$. Then $V(J) \cap \IT^m$ is smooth. (If $a = (\uzero, \uone)$ is
the distinguished point of $X$, then $V(J) \subset \IT^m$.) Let 
$\wI_b$ denote the ideal generated by $I$ in the completion
$\wcO_{\IA^n,b} \cong \IF_b \lbr X,Y \rbr$ at any (closed) point
$b$ (notation of \S3.1). If follows from Theorem 6.5(3) that, at any
$b \in V(J) \cap \IT^m$, the vertices of $\frN(\wI_b)$ are given by
$(\al^i, 0) \in \IN^{n-m} \times \IN^m$, $i = 1,\ldots,s$, together
with $k$ elements $(0,\g^j) \in \{0\} \times \IN^m$ of order
$|\g^j| = 1$, where $k$ is the codimension of $V(J)$ in $\IA^m$
at $b$.
\end{remark}

\section{Samuel stratification of a toric variety}
Subsections 7.2 and 7.3 below isolate the properties of the
Hilbert-Samuel function that play an important part in embedded
resolution of singularities (here in the context of a toric or
binomial variety).
Theorem 7.1 asserts that the maximal Samuel stratum of an affine
binomial variety is the simultaneous
equimultiple locus of a standard basis of the binomial ideal, and 
Theorem 7.2 describes the behaviour of the Hilbert-Samuel function
on blowing up with permissible invariant centre. Theorems 7.1, 7.2
are analogues for binomial varieties over perfect fields of theorems in 
\cite{Hid}, \cite{BMjams}, \cite{BMinv}, \cite{BMsam}.

Subsections 7.4 and 7.5 describe the structure of the Samuel 
stratification of a general embedded toric variety $X$ over a perfect
field. In particular, we prove Theorem 1.6 and we show that the
components of the maximal Samuel stratum satisfy the conditions for
the centres of blowing up given in our main theorem 1.1.

\subsection{Hilbert-Samuel function} Let $R$ denote a Noetherian
local ring with maximal ideal $\um$. The {\it Hilbert-Samuel function}
$H_R: \IN \to \IN$ is defined by
$$
H_R(l)\ =\ \length \frac{R}{\um^{l+1}}, \qquad l \in \IN\ .
$$
If $R$ is a $\uk$-algebra, then $H_R(l) = \dim_{\uk} R/\um^{l+1}$,
for all $l$. We partially order the set of functions $\IN^\IN :=
\{H: \IN \to \IN\}$ as follows: If $H, H' \in \IN^\IN$, then $H \leq H'$
means that $H(l) \leq H'(l)$, for all $l \in \IN$. 

Let $I$ denote an ideal in a formal power series ring
$\uk\lbr x \rbr = \uk\lbr x_1, \ldots, x_n\rbr$. Then, for all
$l \in \IN$, 
\begin{equation}
H_{\uk\lbr x \rbr / I}(l)\ =\ \#\{\al \in \IN^n:\ \al \notin \frN(I),
|\al| \leq l\}
\end{equation}
(by Hironaka's formal division theorem; cf. \S6.3 and Theorem 6.5).
It follows from (7.1) that (for fixed $n$), any non-increasing 
sequence of Hilbert-Samuel functions $H_{\uk\lbr x \rbr / I}$
stabilizes \cite[Thm. 5.2.1]{BMjams}.

If $I \subset \uk\lbr x \rbr$ is a principal ideal of order $\mu$,
then, by (7.1), $H_{\uk\lbr x \rbr / I}(l) = \binom{n+l}{n}$, if
$l < \mu$, and $H_{\uk\lbr x \rbr / I}(l) = \binom{n+l}{n} -
\binom{n+l-\mu}{n}$, if $l \geq \mu$.

We define the Hilbert-Samuel function $H_{X,b}$ of a Noetherian 
local-ringed
space $X = (|X|, \cO_X)$ at a point $b$ as the Hilbert-Samuel
function of the local ring $\cO_{X,b}$.

\subsection{Hilbert-Samuel function and equimultiple locus}
Let $I \subset \uk[x,y] 
= \uk[x_1, \ldots, x_{n-m}, y_1, \ldots, y_m]$ denote a binomial ideal,
where $a = (\uzero, \uone)$ is the distinguished point of the 
corresponding affine binomial variety $X \subset \IA^n$. 

Let $J \subset \uk[y^\pm]$ denote the ideal generated by the binomials
in $I$ involving $y$ alone, and
let $Q$ denote the smooth binomial subvariety of $\IA^n$ defined by the
ideal $J\cdot \uk[x,y^\pm]$. Let
$$
F^i\ =\ x^{\al^i} - x^{\be^i}y^{\g^i}\ \in\ I, \quad i=1, \ldots, s\ ,
$$
denote the standard basis of $I \mod J\cdot \uk[x,y^\pm]$, as in 
Theorem 6.5. For each $i$, if $|\al^i| = 1$, then $x^{\al^i}
= x_{j(i)}$, for some $j(i) \in \{1, \ldots, n-m\}$, and $x_{j(i)}$
occurs (to nonzero power) in no monomial $x^{\al^j}$, $j \neq i$, and
in no monomial $x^{\be^j}$. 

Let $N = N(I)$ 
denote the smooth binomial subvariety of $Q$ defined by the binomials
$F^i(x,y)$, for all $i$ such that $|\al^i| = 1$; i.e., $N$ is the smooth
binomial subvariety of $\IA^n$ defined by the binomials $1 - y^\g \in J$ 
together with $x^{\al^i} - x^{\be^i}y^{\g^i}$, for all $i=1, \ldots, s$
such that $|\al^i| = 1$.

Then $N$ is a minimal embedding submanifold of $X$ (i.e, a smooth
variety of smallest dimension in which $X$ can be embedded; $\dim N$
is determined by $H_{X,a}$).

After reordering the indices $i$ if necessary, we can assume that
\begin{enumerate}
\item
$|\al^i| \geq 2$, $i = 1, \ldots, t$, and $|\al^i| = 1$, 
$i = t+1, \ldots, s$, where $t \leq s$; 
\item
$x_1, \ldots, x_r$,
where $r \leq (n-m)-(s-t)$, are those variables occurring (to nonzero
power) in some $x^{\al^i}$, $i = 1, \ldots, t$.
\end{enumerate}

The mapping $b \mapsto H_{X,b}$ from $X$ to $\IN^\IN$ is Zariski
upper-semicontinuous (cf. \cite{Be}, \cite{BMjams}, \cite{BMinv}).
We define the {\it Samuel stratum} $S_X(b)$ of $b \in X$ as $S_X(b)
:= \{c \in X: H_{X,c} = H_{X,b}\}$. 

\begin{theorem}
\begin{enumerate}
\item
$H_{X,b} \leq H_{X,a}$, for all $b \in X$.
\item
The Samuel stratum $S_X(a)$ is the closed subset $S_{\{F^i\}}(a)$ of
$Q$ given by the simultaneous equimultiple locus of the binomials
$F^i(x,y) = x^{\al^i} - x^{\be^i}y^{\g^i}$, $i=1, \ldots, s$; i.e.,
the closed subset of $N$ defined by
\begin{alignat*}{2}
x_j &= 0,                  & \quad j &= 1, \ldots, r,\\
\mu_{\cdot}(x^{\be^i}) &\geq |\al^i|,  & \quad i &= 1, \ldots, t.
\end{alignat*}
\end{enumerate}
\end{theorem}

\begin{proof}
First suppose that $X$ is toric.
Since every orbit of $X$ is adherent to the distinguished point
$a$, it is enough to prove the assertions in some neighbourhood
of $a$. (1) is then a restatement of semicontinuity. (2) can be
proved exactly in the same way as \cite[Theorem 5.3.1]{BMjams}, using
Theorem 3.1 above (cf. Remark 3.7).

Now consider a general affine binomial variety $X$. As above, (1) and
(2) hold in a neighbourhood of the distinguished point $a$. For every
affine subspace $\IA^q$ of $\IA^n$, $X \cap \IT^q$ is adherent to
$a$, by Lemma 6.3. It therefore suffices to show that the Hilbert-Samuel
function is constant on $X \cap \IT^q$.

To prove the latter, we use the notation of the proof of Lemma 6.3. Let
$K \subset \uk[v^\pm, y^\pm]$ denote the ideal generated by the binomials
in $I$ which involve $(v,y)$ alone. As in \S6.5, let 
$B = \uk[v^\pm, y^\pm]/K$ and let $I_K$ denote the ideal generated by
$I$ in $B[x] \subset B \lbr x \rbr$. If follows from (7.1) and 
Remark 6.7 that, at any point of $X \cap \IT^q$, the Hilbert-Samuel
function of $X$ is completely determined by the codimension of 
$X \cap \IT^q$ in $\IA^q$ and the vertices of $\frN(I_K) \subset
\IN^{n-q}$.
\end{proof}

\begin{corollary}
The Samuel stratum $S_X(a)$ has only simple normal crossings.
Each component of $S_X(a)$ is the intersection with $N$ of an affine
subspace of $\IA^n$ that is transverse to $N$. In particular, in the
toric case, each component of $S_X(a)$ is 
the closure of an orbit of $T_N$
(acting on $N$).
\end{corollary}

\subsection{Hilbert-Samuel function and strict transform}
We continue to use the notation of the preceding subsection.
Let $C$ denote a component of $S_X(a)$. The $C = D \cap N$,
where $D = \{x_j = 0: j \in \D\}$, $\{1, \ldots, r\} \subset \D
\subset \{1, \ldots, n-m\}$, and $\D$ includes no $x_j(i) = 
x^{\al^i}$ with $|\al^i| = 1$ (i.e., with $i = t+1, \dots, s$).

Let $\pi: U' \to \IA^n$ denote the blowing-up of $\IA^n$ with
centre $D$. Let $N' = N(I)'$ and $X'$ denote the strict 
transforms of $N = N(I)$ and
$X$ (respectively) by $\pi$. 
$U'$ is covered by affine charts $U_{x_i} \cong \IA^n$,
$i \in \D$, where, for each $i$, $U_{x_i}$ has affine coordinates
$(x_1, \ldots, x_{n-m}, y_1, \ldots, y_m)$ with respect to which
$\pi$ is given by the following substitution rules: For each $j \in
\D \backslash \{i\}$, substitute $x_i x_j$ for $x_j$, and leave the
remaining coordinates $x_j$ (and all y-coordinates) unchanged. 
Clearly, $X' \cap U_{x_i} \subset N' \cap U_{x_i}$ are
affine binomial subvarieties
of $U_{x_i}$, for each $i$. In the toric case, $X' \subset N'$ are toric
subvarieties of $U'$.

Fix $i \in \D$. Let $I'$ denote the ideal of 
$X' \cap U_{x_i}$. Let $a'$ denote the distinguished
point of $X' \cap U_{x_i}$, and let $J' \subset I'$ denote the
analogue of $J$ above. By general properties of the standard basis
(see \cite[Lemma 3.22]{BMinv}), $I'$ is generated by $J'$ together
with the ``strict transforms'' ${F^j}' = x_i^{-|\al^i|} F^j$ of all
binomials $F^j$. Each ${F^j}'$ is a binomial
$$
{F^j}'\ =\ x^{{\al^j}'} - x^{{\be^j}'}y^{{\g^j}'}\ .
$$
(Certain ${F^j}'$ may belong to $J'$.)

\begin{theorem}
In the preceding notation,
\begin{enumerate}
\item
$H_{X',a'} \leq H_{X,a}$.
\item
$H_{X',a'} = H_{X,a}$ if and only if $a' \in N(I)' = N(I')$,
$|{\al^i}'| = \mu_{a'}({F^i}') = |\al^i|$, $i = 1, \ldots, s$,
and the ${F^i}'$, $i = 1, \ldots, t$, are the elements of
the standard basis
of $I' \mod J'$ of order $\geq 2$.
\end{enumerate}
\end{theorem}

The proof follows that of \cite[Theorem 7.3]{BMjams}.

\subsection{Samuel stratification}
Let $M$ be a smooth variety.
Let $\cE$ denote a finite collection of smooth (Zariski-) closed
subsets of $M$ having only normal crossings. Assume that $\cE$
includes $M$, and that $\cE$ is closed under intersection.
(For example, let $M$ be a smooth toric variety, and let $\cE$ denote
the collection of orbit closures of $M$.)

\begin{lemma}[{\rm cf. \cite[Cor. 1.17]{BMda1}}]
Let $Y$ be a closed subset of $M$. Assume that:
\begin{enumerate}
\item
$Y$ has only normal crossings.
\item
Every component $Z$ of the germ $Y_a$ of $Y$ at any point $a \in Y$
is of the form $Z = (Y \cap E)_a$, where $E \in \cE$.
\end{enumerate}
Then each irreducible component of $Y$ is smooth.
\end{lemma}

\begin{proof}
Let $a \in Y$. Write the germ $Y_a$ as a union of irreducible
components $Y_a = \cup Z_E$, where $Z_E = (Y \cap E)_a$ and $E$
is the intersection of all elements of $\cE$ containing $Z_E$.
We will show that each $Z_E$ extends to an irreducible component
of $Y \cap E$.

Consider any total ordering of $\cE$. Define $e(a) := \max \{E:\ 
Z_E$ is a component of $Y_a\}$. Clearly, $a \mapsto e(a)$
is (Zariski-) upper semicontinuous on $Y$, and its maximum locus
is smooth.

On the other hand, given $a$ and a component $Z_E$ of $Y_a$, we
can choose the ordering of $\cE$  so that $E = \max \cE$. It follows
that $Z_E$ extends to a closed smooth subset of $Y$.
\end{proof}

\begin{example}
Let $M$ be a smooth toric variety over a perfect field $\uk$, 
and let $\cE$ denote the collection of orbit closures of $M$.
Let $X$ be a toric subvariety of $M$, and let $S$ denote the 
maximal Samuel stratum of $X$ (i.e., the locus of maximum values
of the mapping $b \mapsto H_{X,b}$). Then, for every
$E \in \cE$, $Y := S \cap E$ satisfies the hypotheses of Lemma
7.4 (by Theorem 7.1).
\end{example}

Theorem 1.6 is an immediate consequence of Theorem 7.1 and
Lemma 7.4.

Theorem 1.6 generalizes Proposition 5.2. Unlike the case of a 
toric hypersurface, however, in general it is not possible to
desingularize $X$ by blowing up with each successive centre an
arbitrary component of the maximal Samuel stratum.

\begin{example}[{\rm J. Adamus}]
Let $X$ denote the toric subvariety of $\IA^6$ defined by the
ideal $I$ generated by the binomials
\begin{align*}
F &= u^d - x^{d-1}y^dz^{d+1}\ ,\\
G &= v^d - x^dy^{d+1}z^{d-1}\ ,\\
H &= w^d - x^{d+1}y^{d-1}z^d\ ,
\end{align*}  
where $d \geq 2$. (In the notation above, $J = 0$ and $F,G,H$ form
the standard basis of $I$.) The maximal Samuel stratum of $X$ has
$3$ components, each given by setting $u,v,w$ to $0$ along with any
$2$ of the remaining variables $x,y,z$.

Consider the strict transform $X_1 = X'$ of $X_0 = X$ by blowing up
$\IA^6$ with centre any of these components; say, 
$u=v=w=y=z=0$. Then $X_1 \subset U_y \cup U_z$, where $U_y, U_z$
are affine charts defined as in \S8.3. (For example, $U_y \to
U := \IA^6$ is given by substituting $x,,y, yz, yu, yv, yw$ 
for $x,,y,z,u,v,w$ (respectively). Now consider two further blowings-up
whose centres are (components of the maximal Samuel strata) given in
$U_y$ by $u=v=w=x=z=0$, and then in $U_{yz} := (U_y)_z$ by
$u=v=w=y=z=0$. Then the toric ideal of the strict transform $X_3$
in the affine chart $U_{yzz}$ is generated again by the original
binomials $F,G,H$.
\end{example}

\subsection{Combinatorial nature of the Samuel strata}
Let $M = M(\Sig)$ denote a smooth toric variety over a perfect 
field $\uk$, corresponding to a regular fan $\Sig$ in a lattice
$\cong \IZ^n$ (cf. \S2.2). We recall that $M$ is covered by 
affine toric varieties $U_\s$ determined by the maximal cones
$\s$ of $\Sig$.

If $\s, \tau$ are maximal cones, then $Z_{\s \tau} := U_\s \backslash
U_{\tau}$ has only normal crossings. $Z_{\s \tau}$ is the union of the
orbit closures of $U_\s$ that are not in $U_{\tau}$. Moreover, 
$Z_{\s \tau} = \cup E_j$, where the $E_j$ are the intersections with
$Z_{\s \tau}$ of the orbit closures of $M$.

\begin{remark}
If $N$ is a smooth toric subvariety of $M$, then:
\begin{enumerate}
\item
Every orbit closure of $N$ is the intersection with $N$ of an
orbit closure of $M$ that is transverse to $N$. 
(This follows from the affine case.)
\item
The intersection with $N$ of an orbit closure of $M$ has only
normal crossings, and each of its irreducible components is an
orbit closure of $N$ (in particular, smooth).
\end{enumerate}
\end{remark}

The following lemma generalizes Remark 7.7.

\begin{lemma}
Let $X$ denote a closed toric subvariety of $M = M(\Sig)$, and
let $S$ denote the maximal Samuel statum of $X$. Let $E$ be an
orbit closure of $M$ and let $C$ be an irreducible component of
$S \cap E$. Then:
\begin{enumerate}
\item
$C$ is smooth.
\item
If $D$ denotes the smallest orbit closure of $M$ containing $C$,
then $X \cap D = C$.
\end{enumerate}
\end{lemma}

\begin{proof}
(1) is true, by Lemma 7.4 and Example 7.5.

Let $D$ be the smallest orbit closure of $M$ containing $C$.
We first show that $C$ is open and closed in $X \cap D$: 
By Theorem 7.1, in any
chart $U_\s$ where the distinguished point $a_\s$ of $X$ belongs
to $C$, $C \cap U_\s = D \cap S \cap U_\s = D \cap N_\s$,
where $N_\s$ is the minimal embedding submanifold defined before
(depending on the ordering of the affine coordinates); therefore,
$C \cap U_\s = D \cap X \cap U_\s$. In other words, if $C \cap U_\s \neq
\emptyset$, then $C \cap U_\s = (X \cap D) \cap U_\s$.

Now consider two charts $U_\s, U_\tau$. If $C \cap (U_\s \cap U_\tau)
\neq \emptyset$, then $C \cap U_\tau \neq \emptyset$, so that
$C \cap U_\tau = (X \cap D) \cap U_\tau$. On the other hand, suppose that 
$C \cap U_\s \subset Z_{\s \tau}$. Then $C \subset E_j$ for some $j$,
so that $D \subset E_j$, by the definition of $D$. Therefore,
$D \cap U_\tau = \emptyset$; i.e., $(X \cap D) \cap U_\tau = \emptyset$.
This proves (2).
\end{proof}

\section{Proof of the main theorem}
In this section, we prove Theorem 1.1. Let
$X \hookrightarrow M$ denote an equivariant embedding of toric
varieties, where $M$ is smooth.
If $H \in \IN^\IN$, we let $S_H(X)$ denote the Samuel stratum
$S_H(X) := \{b \in X:\ H_{X,b} = H\}$. (In particular, if $a \in X$,
then $S_X(a) = S_{H_{X,a}}(X)$.)

Let $H$ denote the maximum Hilbert-Samuel function of $X$. We will
say that a sequence of blowings-up (1.1) is {\it $H$-permissible}
if conditions (1) and (2) of Theorem 1.1 are satisfied and, in 
addition, $C_j \subset S_H(X_j)$, $j = 0, \ldots, t$ (i.e., each
$X_j$ has maximum Hilbert-Samuel function $H$, and each $C_j$ lies
in the maximal Samuel stratum of $X_j$).

\begin{theorem}
Let $H$ denote the maximum Hilbert-Samuel function of $X$. Then
there is an $H$-permissible sequence of blowings-up (1.1) such that
$S_H(X_{t+1}) = \emptyset$.
\end{theorem}

By the stabilization theorem for the Hilbert-Samuel function
(\cite[Thm. 5.2.1]{BMjams}; cf. \S7.1 above), it follows from Theorem
8.1 that there is a finite sequence of blowings-up (1.1) satisfying
conditions (1), (2) of Theorem 1.1, and also:
\begin{enumerate}
\item[($3'$)]
Each $C_j \subset \Sing X_j$.
\item[($4'$)]
$X_{t+1}$ is smooth.
\end{enumerate}

We prove Theorem 8.1 below. Since each (reduced) component of
the exceptional divisor $E_{t+1}$ is a codimension one orbit closure
of $M_{t+1}$, it will then follow from Theorem 8.9 (in \S8.4 below)
that, beginning with a sequence
of blowings-up (1.1) satisfying (1), (2), ($3'$) and ($4'$), we can
make a further sequence of blowings-up where each $C_j \subset
X_j \cap E_j$, after which (4) is satisfied.

\subsection{Marked monomial ideal}
It is convenient to use some of the structure of W{\l}odarczyk \cite{W};
in particular, his notion of ``marked ideal'' in a simple monomial
setting. Note, however, that our {\it marked monomial ideals} below
are marked by more structure than the marked ideals of \cite{W}.

\begin{definition}
A {\it marked monomial ideal} is a quintuple
$$
\ucH = (M,N,P,\cH,e)\ ,
$$
where
\begin{enumerate}
\item[$M$] is a smooth toric variety,
\item[$N$] is a smooth closed toric subvariety of $M$,
\item[$P$] is a smooth (closed) invariant subvariety of $N$,
\item[$\cH$]$= \cH_1 + \cdots + \cH_r$, where each $\cH_i \subset
\cO_M$  is a product of principal ideals defining codimension
one orbit closures of $M$ that simultaneously have only normal crossings 
with $N$ and do not contain $P$,
\item[$e$] is a positive integer.
\end{enumerate}
\end{definition}

A {\it monomial ideal} $\cH \subset \cO_M$ means an ideal of the
form $\cH_1 + \cdots + \cH_r$, where each $\cH_i \subset \cO_M$ is
a product of principal ideals defining codimension one orbit closures
of $M$.

Note that, in Definition 8.2, if $\dim N > 0$, then the
codimension one orbit closures of $M$ that are involved in $\cH$
are transverse to $N$. It follows that the restriction of $\cH$
to $N$, $\cH \cdot \cO_N$ is a monomial ideal in $\cO_N$ and, for
all $a \in P$, $\mu_a(\cH \cdot \cO_P) = \mu_a(\cH \cdot \cO_N)
= \mu_a(\cH)$ (where $\mu_a$ denotes the order at $a$).

Let $\ucH = (M,N,P,\cH,e)$ be a marked monomial ideal. We define
the {\it support} of $\ucH$,
$$
\supp \ucH := \{a \in P:\ \mu_a(\cH) \geq e\}\ .
$$
Then $\supp \ucH$ is a closed subset of $N$ that has only normal 
crossings. (It is a union of orbit closures.) By Lemma 7.4, every
irreducible component of $\supp \ucH$ is smooth; in fact, every
irreducible component of $\supp \ucH$ is an orbit closure of $N$
(cf. \S7.5). 

We say that $\ucH$ {\it has maximal order} if $\mu_a(\cH) \leq e$,
for all $a \in P$.

Let $\pi = \pi_D:\ M' \to M$ be a blowing-up with centre $D$, where $D$
is a smooth invariant subvariety of $M$. 
We say that $\pi$ is {\it permissible} for
$\ucH$ if $C := D \cap N$ is smooth and $C \subset \supp \ucH$.

\begin{definition}
Let $\ucH = (M,N,P,\cH,e)$ be a marked monomial ideal, and let
$\pi = \pi_D:\ M' \to M$ be a blowing-up that is permissible for 
$\ucH$. The {\it transform} $\ucH'$ of $\ucH$ by $\pi$ is a marked
monomial ideal $\ucH' = (M',N',P',\cH',e')$, where
\begin{enumerate}
\item[$N'$] is the strict transform of $N$ by $\pi$ (so that
$\pi|N':\ N' \to N$ is the blowing-up $\pi_C$ of $N$ with centre
$C = D \cap N$),
\item[$P'$] is the strict transform of $P$,
\item[$\cH'$]$= \cI_{\pi^{-1}(D)}^{-e}\cdot\pi^*(\cH)$, where 
$\cI_{\pi^{-1}(D)} \subset \cO_{M'}$ denotes the principle
ideal defining the exceptional divisor $\pi^{-1}(D) \subset M'$
(a smooth invariant hypersurface in $M'$),
\item[$e'$]$= e$.
\end{enumerate}
\end{definition}

Note that $\cH'|N' = \cH'\cdot\cO_{N'}$ coincides with
$\cI_{(\pi|N')^{-1}(C)}^{-e}\cdot(\pi|N')^*(\cH|N)$.

A {\it permissible} sequence of blowings up for $\ucH$ means a sequence
of blowings-up
\begin{equation}
M = M_0 \stackrel{\pi_1}{\longleftarrow} M_1 \longleftarrow \cdots
\stackrel{\pi_{t+1}}{\longleftarrow} M_{t+1}\ ,
\end{equation}
where, for each $j = 0, \ldots, t$, $\pi_{i+1}$ is permissible
for $\ucH_i = (M_i,N_i,P_i,\cH_i,e_i)$ and 
$\ucH_{i+1} = (M_{i+1},N_{i+1},P_{i+1},\cH_{i+1},e_{i+1})$ denotes
the transform of $\ucH_i$ by $\pi_{i+1}$. (We set $\ucH_0 = \ucH$.)

A {\it resolution of singularities} of $\ucH$ means a sequence of
permissible blowings-up (8.1) such that $\supp \ucH_{t+1} = \emptyset$.

\begin{example}
Suppose that $M$ is affine; say $X \subset M \subset \IA^n$. We use 
the notation of \S7.2 above. (In particular, $N \subset \IA^n$ denotes
a minimal embedding submanifold of $X$.) We can assume that 
$N \subset M$. Consider the standard basis
elements $x^{\al^i} - x^{\be^i}y^{\g^i}$ of $I \mod J\cdot \uk[x,y^\pm]$
of orders $|\al^i| \geq 2$. We can assume that the variables $x$
are listed in two blocks $x = (z,u)$, where $z$ consists of the
{\it essential variables} of the initial monomials $x^{\al^i}$;
i.e., those x-variables which occur (to nonzero power) in some
$x^{\al^i}$, $|\al^i| \geq 2$. So, for each $i$, we write $x^{\al^i}
= z^{\al^i}$ and $x^{\be^i} = z^{\xi^i}u^{\eta^i}$. Set $e_i :=
|\al^i| - |\xi^i|$. 

Each $x_j = 0$ is a codimension one orbit closure of $\IA^n$ that
intersects $M$ (respectively, $N$) transversely in a codimension
one $T_M$ (respectively, $T_N$) orbit closure. Define $P \subset N$
by $z = 0$, set $e := \prod e_i$, and let $\cH \subset \cO_M$ denote
the ideal generated by the monomials $(u^{\eta^i})^{\prod_{j \neq i}e_j}$.
Then $\ucH = (M,N,P,\cH,e)$ is a marked monomial ideal. 

Let $a \in X$ be the distinguished point. Then $H := H_{X,a}$ is the
maximum Hilbert-Samuel function of $X$. It follows from Theorems 7.1
and 7.3 that a sequence of blowings-up (8.1) of $M$ is $H$-permissible
if and only if it is permissible for $\ucH$. Moreover, if (8.1) is
a resolution of singularities of $\ucH$, then $S_H(X_{t+1}) = \emptyset$.
\end{example}

\begin{theorem}
Let $\ucH = (M,N,P,\cH,e)$ be a marked monomial ideal. Then $\ucH$
admits a resolution of singularities. 
\end{theorem}

We will prove this theorem in \S8.4. In \S8.3, we reduce Theorem
8.1 to Theorem 8.5. In these theorems, $M = M(\Sig)$ denotes a smooth
toric variety over a perfect field $\uk$, corresponding to a fan $\Sig$,
and $X$ denotes a closed toric subvariety of $M$.

\begin{remark} Theorem 8.5 in the special case $P = N$ and $e = 1$
is the {\it principalization theorem} for a marked monomial ideal
(cf. \cite{BMda1, V, W}). The principalization theorem for a toric
(or binomial) ideal follows easily from Theorem 8.5 (without using the
Hilbert-Samuel function and Theorems 7.1, 7.3), and implies a weaker
version of the main theorem 1.1 (without condition (3) or normal
flatness), as in \cite{EnV, W}.
\end{remark}

\subsection{Reduction to resolution of singularities of a marked
monomial ideal}

\begin{lemma}
Theorem 8.5 implies Theorem 8.1.
\end{lemma}

\begin{proof}
Assume Theorem 8.5. Let $H$ denote the maximal Hilbert-Samuel function 
of $X$. As in Section 2, we cover $M$ by the affine toric varieties
$U_\s$ corresponding to the maximal cones $\s \in \Sig$. Set $X_\s
:= X|U_\s$. For each maximal cone $\s$, let $\ucH_\s = (U_\s, N_\s,
P_\s, \cH_\s, e_\s)$ denote the marked monomial ideal as defined in
Example 8.4 for $X_\s \subset U_\s$. Number the (finitely many) maximal
cones: $\s^{(1)}, \s^{(2)}, \ldots .$

By Theorem 8.5, there exists a resolution of singularities of
$\ucH^{(1)} := \ucH_{\s^{(1)}}$. By Lemmas 7.4, 7.9 (applied recursively),
there is an $H$-permissible sequence (8.1) 
of blowings-up of $M$ such that,
over $U_{\s^{(1)}}$, this sequence restricts to the preceding resolution 
of singularities of $\ucH^{(1)}$.

This sequence of blowings-up of $M$, restricted to $U_{\s^{(2)}}$,
is permissible for $\ucH_{\s^{(2)}}$. Let $\ucH^{(2)}$ denote the
transform of $\ucH_{\s^{(2)}}$ by this sequence.

We can now resolve the singularities of $\ucH^{(2)}$ using
Theorem 8.5. As before, the sequence of blowings-up involved is
obtained by restriction of an $H$-permissible sequence of blowings-up
of $M^{(1)} := M_{t+1}$. Each centre of blowing-up has empty intersection
with the inverse image of $U_{\s^{(1)}}$.

Now let $\ucH^{(3)}$ denote the transform of $\ucH_{\s^{(3)}}$ and
continue in the same way $\ldots$
\end{proof}

\subsection{Resolution of singularities of a marked monomial ideal}

\begin{definition}
{\it Sum of marked monomial ideals.}
Consider marked monomial ideals $\ucH_1 = (M,N,P,\cH_1,e_1)$ and
$\ucH_2 = (M,N,P,\cH_2,e_2)$. We define
\begin{equation*}
\ucH_1 + \ucH_2\, =\, (M,N,P,\, \cH_1^{e_2}+\cH_2^{e_1},\, e_1e_2)\ .
\end{equation*}
\end{definition}

It is easy to check the following lemma.

\begin{lemma}\hfill
\begin{enumerate}
\item $\supp(\ucH_1 + \ucH_2) = \supp \ucH_1 \cap \supp \ucH_2$.   
\item A sequence of blowings-up of $M$ is permissible for 
$\ucH_1 + \ucH_2$ if and only if it is permissible for both $\ucH_1$
and $\ucH_2$.
\item The transforms by such a sequence satisfy $(\ucH_1 + \ucH_2)'
= \ucH_1' + \ucH_2'$.
\end{enumerate}
\end{lemma}

\begin{proof}[Proof of Theorem 8.5]
The proof is by induction on $\dim P$. The assertion is trivial in
the case $\dim P = 0$. 
\smallskip

\noindent
{\it Step 1.} {\it Resolution of singularities of a marked monomial ideal}
$\ucH = (M,N,P, \cH, e)$ {\it of maximal order.} Let $\s$ denote a maximal
cone of $\Sig$, and set $\ucH_\s = \ucH|U_\s := (U_\s, N_\s, P_\s, \cH_\s,
e)$, where $N_\s := N \cap U_\s$, $P_\s := P \cap U_\s$, $\cH_\s := \cH|U_\s$.
Then $\cH_\s$ is generated by monomials in $u$, where $u$ denotes a certain
block of the affine coordinates which simultaneously have only normal crossings
with $N_\s$, and none of which vanish on $P_\s$ (cf. Example 8.4). 
By Remark 7.8 or Lemma 7.9
(with $X=N$), any permissible sequence of blowings-up for $\ucH_\s$ is the
restriction of a permissible sequence of blowings-up for $\ucH$.

Of course, $\supp \ucH \cap U_\s = \emptyset$ unless $\cH_\s$ is generated 
by monomials of degree $\geq e$. In the latter case, $a_\s \in \supp \ucH$,
where $a_\s$ denotes the distinguished point of $N_\s$.

Suppose that $\supp \ucH \cap U_\s \neq \emptyset$. Then 
$\mu_{a_\s}(\cH_\s) = e$,
so that $\cH_\s$ is not generated by monomials all of degree $> e$.
Say that $\cH_\s$ is generated by monomials $u^{\xi_i}$, $|\xi_i| \geq e$.
Write $u = (z,w)$, where $z$ consists of the ``essential variables'' --
the variables each of which occur (with nonzero power) in some monomial
$u^{\xi_i}$ with $|\xi_i| = e$.

For each $i$, write $u^{\xi_i} = z^{\eta_i}w^{\zeta_i}$. Define 
$Q_\s \subset P_\s$ by $z=0$ (``maximal contact subspace''), and
let $\ucC(\ucH_\s)$ denote the marked monomial ideal
$$
(U_\s, N_\s, Q_\s, \cC(\ucH_\s), e_\s) = \sum_i (U_\s, N_\s, Q_\s,
(w^{\zeta_i}), e - |\eta_i|)\ .
$$
Then
$\supp \ucH_\s = \supp\ucC(\ucH_\s) \subset Q_\s$, a blowing-up of $U_\s$
is permissible for $\ucH_\s$ if and only if it is permissible for 
$\ucC(\ucH_\s)$, and the transforms 
$\ucH_\s'$ and $\ucC(\ucH_\s)'$ by a permissible blowing-up satisfy
$\ucC(\ucH_\s)' = \ucC(\ucH_\s')$. Therefore, 
a sequence of blowings-up of $U_\s$ is permissible for $\ucH_\s$ if and
only if it is permissible for $\ucC(\ucH_\s)$. 

In particular, after any sequence of permissible blowings-up 
for $\ucH_\s$,
writing $\ucH_\s'$ and $\ucC(\ucH_\s)'$ for the transforms, we conclude
that $\supp \ucH_\s' = \supp\ucC(\ucH_\s)' \subset Q_\s'$, and
$\ucH_\s'$, $\ucC(\ucH_\s)'$ have the same sequences of permissible
blowings-up. Moreover (as above), 
any sequence of permissible blowings-up extends to a sequence of
permissible blowings-up for $\ucH'$.

By induction on dimension, since $\dim Q_\s < \dim P_\s$, there
is a sequence of blowings-up of $M$, permissible for $\ucH$, which
restricts to a resolution of singularities of $\ucH_\s$.

We can now argue as in the proof of Lemma 8.7 to complete the 
proof of Step 1:
Number the maximal cones of $\Sig$:
$\s^{(1)}, \s^{(2)},\ldots .$ As above,
there is a sequence of blowings-up of $M$, permissible for $\ucH^{(1)}
:= \ucH$, which
restricts to a resolution of singularities of 
$\ucH_{\s^{(1)}}$.
Let $\ucH^{(2)}$ denote the transform of $\ucH^{(1)}$ by this sequence;
likewise, $(\ucH_{\s})^{(2)}$ and $\ucC(\ucH_{\s})^{(2)}$, for each $\s$.
Write $\pi^{(1)}:\ M^{(2)} \to M^{(1)}=M$ for the composite of the 
sequence of permissible blowings-up.

Suppose that $\supp \ucH^{(2)} \cap (\pi^{(1)})^{-1}(U_{\s^{(2)}}) \neq
\emptyset$. Then $\ucH^{(2)}$ has maximal order, 
$\supp \ucH \cap U_{\s^{(2)}} \neq \emptyset$, 
$(\ucH_{\s^{(2)}})^{(2)}$ and $\ucC(\ucH_{\s^{(2)}})^{(2)}$ have the same
permissible sequences, and each permissible sequence for 
$(\ucH_{\s^{(2)}})^{(2)}$ extends to a
permissible sequence for $\ucH^{(2)}$ whose centres are disjoint from
the inverse images of $U_{\s^{(1)}}$. Therefore,
there is a further permissible sequence of blowings-up, after which
$\supp \ucH^{(3)}$ is disjoint from the inverse image of 
$U_{\s^{(1)}} \cup U_{\s^{(2)}}$. Etc.
\smallskip

\noindent
{\it Step 2.} {\it Resolution of singularities of a marked monomial ideal}
$\ucH = (M,N,P, \cH, e)${\it , in general.} 
\smallskip

\noindent
{\it (a) Reduction to the monomial case.} We can factor $\cH$ as
$\cH = \cD(\cH)\cdot\cQ(\cH)$, where $\cD(\cH)$ is a product of
principal ideals defining codimension one orbit closures transverse
to $N$ (not containing $P$), and $\cQ(\cH)$ is not divisible by any
such principal ideal (nor, then, by any principle ideal 
defining a codimension one orbit closure). Define
$$
\nu_{\ucH} := \max\{\mu_b(\cQ(\cH)):\ b \in \supp \ucH\}\ ,
$$
and set
\begin{align*}
\ucQ(\ucH) &:= (M,N,P, \cQ(\cH), \nu_{\ucH})\ ,\\
\ucD(\ucH) &:= (M,N,P, \cD(\cH), e - \nu_{\ucH})\ .
\end{align*}

We define the {\it companion ideal} of $\ucH$ as the marked monomial
ideal 
$$
\ucG(\ucH) = 
\begin{cases}
\ucQ(\ucH) + \ucD(\ucH)\ , &\text{if $\nu_{\ucH}<e$\ ,}\\
\ucQ(\ucH)\ , &\text{if $\nu_{\ucH} \geq e$}
\end{cases}
$$
(cf. \cite[4.23]{BMinv}, \cite[\S2.6]{BMda1}, \cite[Def.\,3.0.10]{W}). Then
$\ucG(\ucH)$ is a marked monomial ideal of maximal order.

By the {\it monomial case}, we mean that $\cH = \cD(\cH)$.

By Lemma 8.8,
$$
\supp \ucG(\ucH) = \supp \ucQ(\ucH) \bigcap \supp \ucD(\ucH)
= \supp \ucQ(\ucH) \bigcap \supp \ucH\ ,
$$
any permissible sequence of blowings-up for $\ucG(\ucH)$ is permissible
for $\ucH$, and the transforms $\ucH',\, \ucG(\ucH)',\, \ucQ(\ucH)',\, 
\ucD(\ucH)'$ by such a sequence satisfy
\begin{equation}
\supp \ucG(\ucH)' = \supp \ucQ(\ucH)' \bigcap \supp \ucD(\ucH)'
= \supp \ucQ(\ucH)' \bigcap \supp \ucH'\ .
\end{equation}
Note that if $\ucH$ has maximal order, then $\ucH$ and $\ucG(\ucH)$
have the same permissible sequences and the same supports throughout
such a sequence. 

By Step 1 above, there is a resolution of singularities of $\ucG(\ucH)$.
The sequence of blowings-up involved is permissible for $\ucH$. Let
$\ucH_1$ denote the transform of $\ucH$ by this sequence. It follows
from (8.2) that $\nu_{\ucH_1} < \nu_{\ucH}$.
We can now repeat the procedure: let $\ucH_2$ be
the transform of $\ucH_1$ by a resolution of singularities of 
$\ucG(\ucH_1)$, $\ldots$ The process terminates after a finite number
of steps, when we arrive at a marked monomial ideal $\ucH_m$ with
either $\supp \ucH_m = \emptyset$ or $\nu_{\ucH_m} = 0$. In the 
former case, we have a resolution of singularities of $\ucH$, and
in the latter it is enough to resolve the singularities of 
$\ucD(\ucH_m)$; i.e., we have reduced to the monomial case.
\smallskip

\noindent
{\it (b) Monomial case.} $\ucH = (M,N,P, \cH, e)$, where $\cH = \cD(\cH)$.
Suppose that $C$ is an irreducible component of $\supp \ucH$. Then
$C$ is smooth, and $C = N \cap D$, where $D$ is the smallest orbit closure
of $M$ containing $C$. Clearly, if $\ucH' = (M',N',P', \cH', e)$ is the
transform of $\ucH$ by the blowing-up with centre $D$, then $\cH'
= \cD(\cH')$. We will show that, after finitely many such blowings-up,
the transform of $\ucH$ has empty support.

The proof is the same as that of desingularization of a toric 
hypersurface (with $\cD(\cH)$ playing the part of $x^{\g^+}$ in
\S5.1 above), or that of \cite[Case (b), p. 260]{BMinv},
\cite[S3.3(4)]{BMda1}, \cite[Step 2b, p. 20]{W}: Let $\s$ denote
a maximal cone of $\Sig$, and define $\ucH_\s$ as in Step 1 above.
Then $\cH_\s$ is generated by a monomial $u^\om = u_1^{\om_1} \cdots
u_q^{\om_q}$, where $u = (u_1, \ldots, u_q)$ is a block of the
affine coordinates, as in Step 1. If $C \cap U_\s \neq \emptyset$, then
$C \cap U_\s \subset P_\s$ can be described as $u_{j_1} = \cdots
= u_{j_l} = 0$, where $1 \leq j_k \leq q$ ($k = 1,\ldots,l$),
\begin{align*}
\om_{j_1} + \cdots + \om_{j_l} &\geq e\ ,\\
\om_{j_1} + \cdots + {\hat \om_{j_k}} + \cdots + \om_{j_l} &< e\ ,
\quad k = 1,\ldots,l\ ,
\end{align*}
(where ${\hat \om_{j_k}}$ means that $\om_{j_k}$ is deleted).
After blowing up with centre $D$, 
the strict transform $(N_\s)'$ lies in the union of the
affine charts corresponding to the variables $u_{j_k}$. In the chart
corresponding to $u_{j_k}$, $\cH_\s'$ is generated by $u^{\om'}$,
where
\begin{align*}
\om_j' &= \om_j\ , \quad j \neq j_k\ ,\\
\om_{j_k}' &= \om_{j_1} + \cdots + \om_{j_l} - e < \om_{j_k}\ ;
\end{align*}
i.e., $|\om'| < |\om|$. Therefore, the support of the transform
of $\ucH$ is disjoint from the inverse image of $U_\s$ after at
most $|\om|-e+1$ such blowings-up over $U_\s$.
\end{proof}

\subsection{Transformation to normal crossings} The following
result is needed to complete our proof of Theorem 1.1.

\begin{theorem}
Let $X \hookrightarrow M$ denote an equivariant
embedding of smooth toric varieties over a perfect field $\uk$.
Let $\Th$ denote a set of smooth invariant hypersurfaces in $M$.
Then there is a finite sequence
of blowings-up (1.1) such that: 
\begin{enumerate}
\item 
The centre $D_j$ of each blowing-up $\pi_{j+1}$ is a smooth
invariant subvariety of $M_j$.
\item
Set $X_0 = X$ and $\Th_0 = \Th$. 
For each $j = 0,\ldots,t$, let $X_{j+1}$ denote
the strict transform of $X_j$ by $\pi_{j+1}$, and let $\Th_{j+1}$
denote the set of strict transforms of elements of $\Th_j$. Then, for
each $j = 0,\ldots,t$, $C_j := D_j \cap X_j$ \ is a smooth invariant
subvariety of $X_j$, $D_j$ is the smallest invariant subvariety of 
$M_j$ containing $C_j$, and each component of $C_j$ is the intersection
of $X_j$ and certain components of elements of $\Th_j$. 
\item
$X_{t+1}$ is disjoint from all elements of $\Th_{t+1}$.
\end{enumerate}
\end{theorem}

The proof of Theorem 8.10 parallels the preceding part of this section. 
If $a \in X_j$, let $s(a)$
denote the number of elements of $\Th_j$ containing $a$. Let $s :=
\max_{a \in X}s(a)$. For
any sequence of blowings-up (1.1) satisfying conditions (1) and (2) of
Theorem 8.9, let $S_\Th(X_j) := \{a \in X_j: s(a) =s\}$. Then, for 
each $j$, $S_\Th(X_j)$ has only normal crossings. As in Lemma 7.9,
if $E$ is an orbit closure of $M_j$ and $C$ is an irreducible component
of $S_\Th(X_j) \cap E$, then $C$ is smooth and $X_j \cap D = C$, where
$D$ is the smallest orbit closure of $M$ containing $C$.

A sequence of blowings-up (1.1) satisfying conditions (1) and (2) of
Theorem 8.10 will be called $\Th${\it -permissible}. Theorem 8.10 is
an obvious consequence of the following.

\begin{theorem}
There is a $\Th$-permissible sequence of blowings-up (1.1), such
that $S_\Th(X_{t+1}) = \emptyset$.
\end{theorem}

\begin{proof}
First consider $M$ affine, as in Example 8.4. Since $X$ is smooth,
$N = X$ and each standard basis element $x^{\al^i} - x^{\be^i}y^{\g^i}$
of $I \mod J\cdot \uk[x,y^\pm]$ is of order $1$; i.e., $x^{\al^i} =
x_{j(i)}$ for some $j(i)$. We can assume that the variables are listed
in two blocks $x = (z,u)$, where $z$ consists of the variables $x_{j(i)}$,
and each $x^{\be^i}$ is a monomial in $u$ (that we denote $u^{\be^i}$).

Set $P=N$. Each $H \in \Th$ is defined by $z_k = 0$, for some $k$,
or by $u_l = 0$, for some $l$. Let $\cH \subset \cO_M$ denote the
ideal generated by $u^{\be^k}$ (for all such $k$) and $u_l$ (for all
such $l$). Then $\ucH = (M,N,P, \cH, 1)$ is a marked monomial ideal.

A sequence of blowings-up (8.1) of $M$ is $\Th$-permissible if and
only if it is permissible for $\ucH$. Moreover, if (8.1) is a
resolution of singularities of $\ucH$, then $S_\Th(X_{t+1}) = \emptyset$.

Theorem 8.11 follows from Theorem 8.5 as in \S8.2 above.
\end{proof}

\section{Algorithm for canonical equivariant desingularization}
Our proof of Theorem 1.1 in Section 8 provides an algorithm for 
equivariant embedded desingularization of a toric variety, but allows
certain noncanonical choices for the centres of blowing up. The purpose
of this section is to give an algorithm for canonical equivariant
resolution of singularities, as well as an invariant that tracks
the progress of the algorithm towards desingularization. Both the
algorithm and the invariant can be described in a very simple way.
The proofs in this section, however, depend on techniques of \cite{BMinv}
or \cite{BMda1} that we do not redevelop, even though they can be 
considerably simplified for toric (or binomial)
varieties (over an arbitrary perfect
field). We do not see how to obtain the results here using 
W{\l}odarczyk's approach to canonicity \cite{W}, although we use
his language of marked ideals. 

The algorithm presented here corresponds to Theorem 8.1, and thus
to Theorem 1.1 with condition (4) replaced by ($4'$) $X_{t+1}$ is smooth.
(See Theorem 8.1 ff.) The same algorithm can be run for Theorem 8.10
to get all the conditions of Theorem 1.1.

In \S9.1, we describe an algorithm for choosing a centre of blowing up
locally, that is implicit in the proof of Theorem 8.5. The local algorithm
is canonical modulo an ordering of the codimension one orbit closures of
$M$.

Our proofs of Lemma 8.7 and of Step 1 in Theorem 8.5 use the orbit
structure of the torus action (in particularly, Lemma 7.9) to globalize
the local choice of centre of blowing up. In \S9.2, we show how to
replace this argument by a global algorithm that is canonical modulo
an ordering of the codimension one orbit closures. 

In \S9.3, we give an algorithm for canonical equivariant embedded
desingularization (Addendum 1.2), by a simple variation of the 
preceding construction: The only change is that, in Step 2 of Theorem
8.5, the marked monomial ideal $\ucD(\ucH)$ and thus the companion
ideal $\ucG(\ucH)$ are defined using only those codimension one orbit
closures which arise as (reduced) components of the exceptional divisors.
At any stage of the resolution process, the ``exceptional hypersurfaces''
(the strict transforms of the exceptional divisors of all previous
blowings-up) are
ordered by their ``years of birth''. (See Definitions 9.1.)

Let $M = M(\Sig)$ denote a smooth toric variety over a perfect field 
$\uk$, as in Section 8.

\begin{definitions} 
{\it Ordering of invariant hypersurfaces and exceptional divisors.}
Consider any totally ordered subset $\cE$ of the set of smooth invariant
hypersurfaces in $M$; say $\cE = \{H_1, \ldots, H_q\}$. Associate to any
subset $I$ of $\cE$, a finite sequence $\de(I) = 
(\de_1, \ldots, \de_q)$, where, for each $i = 1,\ldots,q$, $\de_i = 0$
if $H_i \not\in I$, and $\de_i = 1$ if $H_i \in I$. We totally order the
subsets $I$ of $\cE$ using the lexicographic order $\text{lex}$ of the
sequences $\de(I)$; i.e., $I_1 < I_2$ means that $\text{lex}\, \de(I_1)
< \text{lex}\, \de(I_2)$. (For example, $\emptyset \leq I$, for all
$I \subset \cE$.)

Suppose that $\pi:\ M' \to M$ is a blowing-up whose centre $D$ is a smooth
invariant subset
of $M$ ($\codim D \geq 2$). Set $\cE' := \{H_1', \ldots, H_q',
H_{q+1}'\}$, where $H_i'$ is the strict transform of $H_i$, 
$i = 1,\ldots,q$, and $H_{q+1}' = \pi^{-1}(D)$.

For example, take $\cE = \emptyset$ and consider any sequence of 
blowings-up (1.1) whose centres are smooth invariant subvarieties. 
Set $\cE_0 := \cE$ and
$\cE_{j+1} := \cE_j'$, $j = 0,\ldots,t$. Thus each $\cE_j$ is the set of
strict transforms (in $M_j$) of the exceptional divisors
$\pi_i^{-1}(D_{i-1})$, $i = 1,\ldots,j$, ordered by ``year of birth''.
\end{definitions}

\subsection{Local calculation of the centre of blowing up} Consider
the affine case $X \subset M \subset \IA^n$, as in Example 8.4. Let
$\cE$ denote the set of codimension one orbit closures of $M$, with
a fixed total ordering. Let $a \in X$ denote the distinguished point.
Let $\ucH = (M,N,P, \cH, e)$ be the marked monomial ideal defined in
Example 8.4. (In the language of \cite{BMinv}, we would call $\ucH$
a ``presentation'' of the Hilbert-Samuel function of $X$ at $a$.)

We will extract from the proof of Theorem 8.1, an algorithm for choosing
a centre of blowing up $D$, as well as a local ``invariants'' $\inv (a)$,
$\mu(a)$ and $J(a)$. $\inv(a)$ is a finite sequence $(H_{X,a},
\nu_2, \ldots, \nu_{q+1})$, where $\nu_2, \ldots, \nu_q$ are positive
rational numbers, $\nu_{q+1} = 0$ or $\infty$, and $q \leq \dim N$ (the
minimal embedding dimension of $X$); $\mu(a)$ is a rational number defined
in the case that $\nu_{q+1} = 0$; $J(a)$ is a subset of $\cE$.

Recall that a sequence of blowings-up (8.1) is $H$-permissible if and
only if it is permissible for $\ucH$. Moreover, if (8.1) is a resolution
of singularities of $\ucH$, then $S_H(X_{t+1}) = \emptyset$.

Therefore, we consider resolution of singularities of an arbitrary
marked monomial ideal $\ucH = (M,N,P, \cH, e)$ in the affine case
$M \subset \IA^n$. We will define a centre of blowing up $D_{\ucH}$, as
well as ``invariants'' $\inv_{\ucH}(a)$, $\mu_{\ucH}(a)$, $J_{\ucH}(a)$,
by induction on $\dim P$.

Then, to define $\inv(a)$ and $D$ for $X \subset M$, we take $\ucH$ 
as in Example 8.4 and set $D := D_{\ucH}$,
$$
\inv(a) := (H_{X,a}, 1, \ldots, 1, \inv_{\ucH}(a))\ ,
$$
where there are $\codim(P \subset N) - 1$ ones, $\mu(a) := \mu_{\ucH}(a)$,
and $J(a) := J_{\ucH}(a)$. ($\codim(P \subset N)$ means the codimension
of $P$ in $N$.)

Let $\ucH = (M,N,P, \cH, e)$ denote a marked monomial ideal, where
$M \subset \IA^n$.
\smallskip

\noindent
{\it (1) Monomial case} $\cH = \cD(\cH)$. Define $\inv_{\ucH}(a) := 0$ 
and $\mu_{\ucH}(a) = \mu_a(\cH)/e$. By Step 2(b) in \S8.3, 
each irreducible
component $Z$ of $\supp \ucH$ is of the form  $Z = N \cap D$, where
$D := \{H_i \in \cE:\ Z \subset H_i\}$ (the smallest orbit closure of 
$M$ containing $Z$). Write $Z = Z_I$ , where 
$I := \{H_i \in \cE:\ Z \subset H_i\}$.

Consider a blowing up $\pi:\ M' \to M$ with centre $D = D_I$ corresponding
to a component $Z_I$ of $\supp \ucH$. Let $a'$ denote a distinguished
point of the strict transform $N'$ (in an affine chart of $M'$, as 
described in Step 2(b)). If $a' \in \supp \ucH'$, then $1 \leq 
\mu_{\ucH'}(a') \leq \mu_{\ucH}(a) - 1/e$.

Take $J_{\ucH}(a) := \max \{I:\ Z_I$ is a component of $\supp \ucH \}$,
and set $D_{\ucH} := D_{J_{\ucH}(a)}$.
\smallskip

\noindent
{\it (2) If $\cH \neq \cD(\cH)$, we consider two cases:}
\smallskip

\noindent
{\it (a) $\ucH$ is of maximal order.} (Then $\cH \neq 0$.) 
Define $\ucC(\ucH) = (M,N,Q, \cC(\cH), e_{\ucC(\ucH)})$ as in 
Step 1 of the proof of Theorem 8.5. (Then $\ucC(\ucH)$ is not of
maximal order.) A resolution of singularities of $\ucC(\ucH)$
is a resolution of singularities of $\ucH$.

Since $\dim Q < \dim P$, $\inv_{\ucC(\ucH)}(a)$, $\mu_{\ucC(\ucH)}(a)$,
$J_{\ucC(\ucH)}(a)$, and a centre of blowing up $D_{\ucC(\ucH)}$ are 
determined by induction. Set $D_{\ucH} := D_{\ucC(\ucH)}$, 
$\inv_{\ucH}(a) := (1, \ldots, 1, \inv_{\ucC(\ucH)}(a))$ (where there
are $\codim (Q\subset P)-1$ ones), $\mu_{\ucH}(a) := \mu_{\ucC(\ucH)}(a)$
(if $\inv_{\ucC(\ucH)}(a)$ ends with $0$),
and $J_{\ucH}(a) := J_{\ucC(\ucH)}(a)$.
\smallskip

\noindent
{\it (b) $\ucH$ is not of maximal order.} If $\cH = 0$ (zero ideal),
set $\inv_{\ucH}(a) := \infty$, $J_{\ucH}(a) := \emptyset$, and
take $D_{\ucH} :=$ the smallest orbit closure $D$ of $M$ such that
$P = N \cap D$. In this case, if $\pi:\ M' \to M$ is the blowing-up
with centre $D$, then $\supp \ucH' = \emptyset$.

On the other hand, if $\cH \neq  0$, then we construct the companion
ideal $\ucG(\ucH)$ as in Step 2(a) of the proof of Theorem 8.5. Then
$\ucG(\ucH)$ is of maximal order. Set 
$$
\inv_{\ucH}(a)\, :=\, \left( \frac{\nu_{\ucH}}{e}\, ,\, \inv_{\ucG(\ucH)}(a)
\right)\ ,
$$
$\mu_{\ucH}(a) := \mu_{\ucG(\ucH)}(a)$ (if $\inv_{\ucG(\ucH)}(a)$ 
ends with $0$),
$J_{\ucH}(a) := J_{\ucG(\ucH)}(a)$, and take 
$D_{\ucH} := D_{\ucG(\ucH)}$.
\smallskip

This completes the inductive definition of $\inv(a)$, $\mu(a)$,
$J(a)$, and the centre of blowing up $D$. After any sequence of
blowings-up (1.1) satisfying conditions (1) and (2) of Theorem 1.1,
we can repeat the algorithm above in any affine chart of $M_{t+1}$, 
using the induced order of $\cE_{t+1}$, where $\cE_0 = \cE$
(Definitions 9.1).

The invariants $\inv$, $\mu$ and $J$ and the centre of blowing up
defined above depend only on the local equivariant isomorphism class 
of $X \subset M$ and the ordering of $\cE$ \cite{BMinv}, 
\cite{BMda1}. The following theorem is easy to prove using the
construction above. (See, for example, \cite[p.260]{BMinv}, 
\cite[\S3.3]{BMda1}.)

\begin{theorem}
Let $a \in X \subset M \subset \IA^n$, and let $\pi:\ 
M' \to M$ denote the blowing-up with centre $D$ (as above). For every
distinguished point $a'$ of $X' \subset M'$, either $\inv(a') < \inv(a)$,
or $\inv(a)$ ends with $0$, $\inv(a') = \inv(a)$ and $\mu(a') < \mu(a)$.
\end{theorem}

\begin{remarks}
The invariant $\inv(a)$ is a simplification of the general 
desingularization invariant $\inv_X(a)$ of \cite{BMinv}, \cite{BMda1}.
The latter is a finite sequence $(H_{X,a}, s_1(a), \nu_2(a), s_2(a),
\ldots, s_q(a), \nu_{q+1}(a))$, where the terms $s_i(a)$ count
certain exceptional hypersurfaces that do not necessarily
have normal crossings with respect to the maximal contact subspaces.
In the case of toric varieties, the normal crossings condition is
automatic, and the terms $s_i(a)$ (reflecting blowings-up needed
to obtain normal crossings) are unnecessary. See
\cite[\S1]{BMda1} for a discussion of this point. Example 1.3 
illustrates the effect on the complexity of the desingularization
algorithm. 

The construction of $\ucC(\ucH)$ (Step 1 in the proof of Theorem 8.5,
and Case 2(a) above) allows passage to a maximal contact subspace $Q$
of arbitrary codimension in $P$, whereas the general construction
of \cite{BMinv}, \cite{BMda1}, \cite{W} involves inductive steps of
codimension $1$ -- the extra $\nu_i$ terms thus involved in the
general definition of $\inv_X$ are the sequences of ones appearing
in the definition of $\inv$ above.
\end{remarks}

\subsection{Global algorithm}
Let $\cE$ denote the set of codimension one orbit closures of $M
= M(\Sig)$, with any fixed total ordering. Let $X \subset M$ denote
a toric subvariety of $M$. We cover $M$ by the affine toric varieties
$U_\s$ corresponding to the maximal cones $\s$ of $\Sig$. For each
maximal cone $\s$, set $X_\s := X|U_\s$ and let $a_\s$ denote the
distinguished point of $X_\s$.

For each $\s$, we consider the invariants $\inv(a_\s)$, $J(a_\s)$,
and the centre of blowing up $D_\s$ constructed as in \S9.1 for
$X_\s \subset U_\s$.

\begin{theorem}
There is a unique smooth closed invariant subspace $D$ of $M$ such
that $D|U_\s = D_\s$ for every maximal cone $\s \in \Sig$ which realizes
the maximum (lexicographic) value of $(\inv(a_\s), J(a_\s))$ (over
the maximal cones $\s$).
\end{theorem}

The invariant subspace $D$ of the theorem is the smallest closed
invariant subspace of $M$ containing the centre of blowing up
that is prescribed by the desingularization algorithm of \cite{BMinv},
\cite{BMda1} in the simplified version sufficient for toric varieties.
Locally, the algorithm works as in \S9.1 above; semicontinuity
properties of the local invariant guarantee that its maximum locus
is a global smooth closed subspace.
\smallskip

\noindent
{\it Desingularization algorithm.} Set $M_0 = M$, $X_0 = X$, and
$\cE_0 = \cE$. Define a sequence of blowings-up (1.1) by applying
Theorem 9.4  successively to each strict transform $X_j \subset M_j$
(with the set $\cE_j$ 
ordered as in Definitions 9.1), in order to define the centre $D_j$
of the next blowing-up $\pi_{j+1}$. 
\smallskip

The algorithm terminates when $X_{t+1}$ is smooth. (It terminates
because of Theorem 9.2.) In order to obtain the normal crossings
condition (4) of Theorem 1.1, we rerun the preceding algorithm in
the context of Theorem 8.10 (rather than Theorem 8.1, as above).

\subsection{Canonical equivariant embedded desingularization}
The desingularization algorithm of the preceding subsection is
canonical modulo the ordering of the set $\cE$ of codimension one
orbit closures of $M$. We can obtain Addendum 1.2 by a simple
change in the algorithm: Take $\cE_0 = \cE := \emptyset$, and
define $\cE_j$ over any sequence of blowings-up (1.1) satisfying
conditions (1) and (2) of Theorem 1.1, as in Definitions 9.1;
i.e., each $\cE_j$ is the set of strict transforms in $M_j$
of the exceptional divisors $\pi_i^{-1}(D_{i-1})$, $i = 1,\ldots,j$, 
ordered by year of birth.

We run the algorithm exactly as before, with one change: At any
stage $j$ of the process, $\cD(\cH)$ is a product of principal ideals
defining codimension one orbit closures (transverse to $N$ and not
containing $P$, in the notation of Step 2(a) of the proof of
Theorem 8.5) {\it which are irreducible components of the elements
of the exceptional set}\, $\cE_j$. For example,
in ``year zero'', $\cE_0 = \emptyset$, so that, if $\ucH = (M,N,P,
\cH, e)$, then $\ucG(\ucH) = \ucQ(\ucH) = (M,N,P, \cH, \nu_{\ucH})$.
(See Theorem 8.5, Step 2(a).) 

All the constructions and proofs in Sections 8 and 9 (above) carry
over with this change. The effect of the change is that extra 
blowings-up are introduced to ``relabel'' orbit closures as exceptional
divisors before they can be factored out to define the companion ideals.
There is an important effect on the complexity of the desingularization
algorithm -- see Example 1.3.

\section{Toroidal and binomial varieties}
\begin{definitions}
Let $X \hookrightarrow M$ be an embedding of algebraic varieties
over $\uk$, where $M$ is smooth. We will say that $X \hookrightarrow M$ 
is:
\begin{enumerate}
\item {\it locally toric} if, for every (closed) point $a \in X$, there
is an open neighbourhood $U$ of $a$ in $M$, an (equivariant) embedding
of affine toric varieties $Y \hookrightarrow V \subset \IA^n_{\uk}$
(where $V$ is smooth), and an \'etale morphism $\eta: U \to V$ such
that $\eta(X|U) = Y|\eta(U)$.
\item {\it toroidal} if there is a collection $\cE = \{H_i\}$ of smooth
irreducible hypersurfaces in $M$ having only normal crossings, and (1)
holds with the additional condition that $\eta^{-1}(T_V) 
= U \backslash \cup H_i$.
\item {\it locally binomial} if we weaken (1) by assuming only that
$Y \hookrightarrow V \subset \IA^n_{\uk}$ is an embedding of affine
binomial varieties (\S6.1).
\item {\it binomial} if there is a collection $\cE = \{H_i\}$ of smooth
irreducible hypersurfaces in $M$ having only normal crossings, and (3)
holds with the additional condition that $\eta^{-1}(V \cap \IT^n)
= U \backslash \cup H_i$.
\end{enumerate}
In each case, the \'etale morphism $\eta$ involved will be called
a {\it local model}.
\end{definitions}

Let $X \hookrightarrow M$ be a toroidal (or binomial) embedding, and
let $\cE = \{H_i\}$ denoted the associated collection of divisors.

\begin{definition}
A smooth subvariety $D$ of $M$ will be called {\it combinatorial} if
\begin{enumerate}
\item each irreducible component of $D$ is an irreducible component
of an intersection of divisors $H_i$;
\item for every local model $\eta: U \to V$, there is a smooth invariant
subset $D_{\eta}$ of $V$ such that $D \cap U = \eta^{-1}(D_{\eta})$.
\end{enumerate}
(In analogy with the toric case, we use ``invariant subset'' of
an affine binomial variety $V \subset \IA^n$ to mean an intersection
of $V$ with a $\IT^n$-invariant subset of $\IA^n$.)
\end{definition}

Consider a blowing-up $\pi: M' \to M$ with combinatorial centre
$D$. Let $X'$ denote the strict transform of $X$ by $\pi$. Let $\cE'$
denote the collection of smooth irreducible hypersurfaces in $M'$
comprising the strict transforms $H_i'$ of all $H_i \in \cE$, together
with the irreducible components of $\pi^{-1}(D)$. Then $X' \hookrightarrow
M'$ is a toroidal (or binomial) embedding (with respect to $\cE'$).

Our analysis of the Samual stratification in Section 7, and proofs of
resolution of singularities in Sections 8, 9 go through for a toroidal
(or binomial) embedding over a perfect field $\uk$, exactly as before,
with components of intersections of divisors in $\cE$ playing the part
of orbit closures in $M$ in the toric case. In particular, we have the
following analogue of Theorem 1.1 and Addendum 1.2.

\begin{theorem}
Let $X \hookrightarrow M$ (where $M$ is smooth) denote a toroidal 
(or binomial) embedding over a perfect field $\uk$. Then there is 
a finite sequence of blowings-up of $M$,
\begin{equation}
M = M_0 \stackrel{\pi_1}{\longleftarrow} M_1 \longleftarrow \cdots
\stackrel{\pi_{t+1}}{\longleftarrow} M_{t+1}\ ,
\end{equation}
such that:
\begin{enumerate}
\item The centre $D_j$ of each blowing-up $\pi_{j+1}$ is combinatorial.
\item Denote by $X_j$ the successive strict transforms of $X_0 = X$.
Then, for each $j=0,\ldots,t$, $C_j := D_j \cap X_j$ is smooth and
$X_j$ is normally flat along $C_j$.
\item For each $j=0,\ldots,t$, either $C_j \subset \Sing X_j$ or
$X_j$ is smooth and $C_j \subset X_j \cap E_j$, where $E_j$ denotes
the exceptional divisor of $\pi_1\circ\cdots\circ \pi_j$.
\item $X_{t+1}$ is smooth and $X_{t+1}, E_{t+1}$ simultaneously have
only normal crossings.
\end{enumerate}

Moreover, for every embedding $X \hookrightarrow M$ as above
(with collection of divisors $\cE$), there
is a sequence of blowings-up (10.1) satisfying conditions (1)-(4),
with the following property: If $\io: M' \hookrightarrow M$ is an
open toroidal (or binomial) embedding (with respect to $\cE$), then
$X$ and $X' = \io^{-1}(X)$ have the same resolution towers over $M'$
(not counting isomorphisms in the sequences of blowings-up).
\end{theorem}

In the case of a locally toric (or locally binomial) embedding over
a perfect field, the general desingularization algorithm of \cite{BMinv,
BMda1} applies to give canonical embedded resolution of singularities,
with the additional property that, for every local model $\eta: U \to V$,
the sequence of strict transforms of $X|U$ is induced by an embedded
toric (or binomial) desingularization of $Y \hookrightarrow V$
(notation of Definitions 10.1). The general algorithm involves additional
complexity. (See Section 1 and Remarks 9.3.)

It would be interesting to find good generalizations of the results
in this paper to the larger class of ``binomial varieties'' of 
\cite{ES} (varieties defined locally by ideals that are generated by
polynomials with at most two terms).

\bibliographystyle{ams}

\end{document}